\documentclass[english]{article}
\usepackage{natbib}
 \bibpunct[, ]{(}{)}{,}{a}{}{,}%

\usepackage{algorithmicx,algorithm}
\usepackage[noend]{algpseudocode}
\usepackage{graphicx}
\usepackage{amsmath}
\usepackage{multirow}
\usepackage{subfigure}
\usepackage{amsmath}
\usepackage{amssymb}
\usepackage{setspace}
\usepackage[T1]{fontenc}
\usepackage[latin9]{inputenc}
\usepackage{geometry}
\usepackage{babel}
\usepackage{mathrsfs}
\usepackage{mathtools}
\usepackage{amsmath}
\usepackage{amsthm}
\usepackage{amssymb}
\usepackage[unicode=true,pdfusetitle,
bookmarks=true,bookmarksnumbered=false,bookmarksopen=false,
breaklinks=false,pdfborder={0 0 0},pdfborderstyle={},backref=false,colorlinks=false]
{hyperref}
\usepackage{breakurl}

\makeatletter
\theoremstyle{plain}
\newtheorem{thm}{\protect\theoremname}[section]
\theoremstyle{definition}

\theoremstyle{plain}
\newtheorem{assumption}[thm]{\protect\assumptionname}
\theoremstyle{plain}

\theoremstyle{remark}

\theoremstyle{plain}

\theoremstyle{plain}
\newtheorem{lem}[thm]{\protect\lemmaname}

\makeatother

\providecommand{\assumptionname}{Assumption}
\providecommand{\corollaryname}{Corollary}
\providecommand{\definitionname}{Definition}
\providecommand{\lemmaname}{Lemma}
\providecommand{\propositionname}{Proposition}
\providecommand{\remarkname}{Remark}
\providecommand{\theoremname}{Theorem}

\begin{document}
%%%%%%%%%%%%%%%%

%\setcounter{page}{1}
%\DOI{10.1287/xxxx.0000.0000}%




% \TITLE{Bundling Information Goods of Decreasing Value}
\title{A Multi-Level Simulation Optimization Approach for Quantile Functions}

\author{Songhao Wang \thanks{Department of Industrial Systems Engineering and Management, National University of Singapore. Email: {\tt wangsonghao@u.nus.edu}}\and Szu Hui Ng \thanks{Department of Industrial Systems Engineering and Management, National University of Singapore. Email: {\tt isensh@nus.edu.sg}} \and William B. Haskell \thanks{Department of Industrial Systems Engineering and Management, National University of Singapore. Email: {\tt wbhaskell@gmail.com}}}

\maketitle

\begin{abstract}
Quantile is a popular performance measure for a stochastic system to evaluate its variability and risk. To reduce the risk, selecting the actions that minimize the tail quantiles of some loss distributions is typically of interest for decision makers. When the loss distribution is observed via simulations, evaluating and optimizing its quantile functions can be challenging, especially when the simulations are expensive, as it may cost a large number of simulation runs to obtain accurate quantile estimators. In this work, we propose a multi-level metamodel (co-kriging) based algorithm to optimize quantile functions more efficiently. Utilizing non-decreasing properties of quantile functions, we first search on cheaper and informative lower quantiles which are more accurate and easier to optimize. The quantile level iteratively increases to the objective level while the search has a focus on the possible promising regions identified by the previous levels. This enables us to leverage the accurate information from the lower quantiles to find the optimums faster and improve algorithm efficiency. 
\end{abstract} % Enter your abstract
\global\long\def\EE{\mathbb{E}}

\global\long\def\PP{\mathbb{P}}

\global\long\def\RR{\mathbb{R}}

\global\long\def\ZZ{\mathbb{Z}}

\global\long\def\NN{\mathbb{N}}

\global\long\def\TT{\mathbb{T}}
\global\long\def\UU{\mathbb{U}}
\global\long\def\VV{\mathbb{V}}
\global\long\def\WW{\mathbb{W}}
\global\long\def\SS{\mathbb{S}}

\global\long\def\FF{\mathbb{F}}
\global\long\def\XX{\mathbb{X}}
\global\long\def\YY{\mathbb{Y}}
\global\long\def\ZZ{\mathbb{Z}}
\global\long\def\AA{\mathbb{A}}

\global\long\def\Ff{\mathcal{F}}
\global\long\def\Hh{\mathcal{H}}
\global\long\def\Vv{\mathcal{V}}
\global\long\def\Cc{\mathcal{C}}
\global\long\def\Ii{\mathcal{I}}

\global\long\def\Ll{\mathcal{L}}
\global\long\def\Pp{\mathcal{P}}
\global\long\def\Oo{\mathcal{O}}
\global\long\def\Mm{\mathcal{M}}

\global\long\def\pP{\mathscr{P}}
\global\long\def\lL{\mathscr{L}}
\global\long\def\bB{\mathscr{B}}
\global\long\def\dD{\mathscr{D}}
\global\long\def\rR{\mathscr{R}}

\global\long\def\RrR{\mathfrak{R}}
\global\long\def\PpP{\mathfrak{P}}
\global\long\def\XxX{\mathfrak{X}}
\global\long\def\YyY{\mathfrak{Y}}
\global\long\def\ZzZ{\mathfrak{Z}}
\global\long\def\ddD{\mathfrak{d}}
\global\long\def\mmM{\mathfrak{m}}
\global\long\def\nnN{\mathfrak{n}}
\global\long\def\TtT{\mathfrak{T}}
\global\long\def\CcC{\mathfrak{C}}
\global\long\def\EeE{\mathfrak{E}}
\global\long\def\AaA{\mathfrak{A}}
\global\long\def\BbB{\mathfrak{B}}

\global\long\def\iiI{\mathfrak{i}}
\global\long\def\jjJ{\mathfrak{j}}
\global\long\def\kkK{\mathfrak{k}}
\global\long\def\llL{\mathfrak{l}}
\global\long\def\mmM{\mathfrak{m}}
\global\long\def\nnN{\mathfrak{n}}

\global\long\def\D{\mathrm{d}}

\global\long\def\E{\mathrm{e}}

\global\long\def\Reals{\RR}
\global\long\def\RealsNN{\RR_{\geq0}}
\global\long\def\RealsP{\RR_{>0}}

\global\long\def\deq{\coloneqq}
\global\long\def\eqd{\eqqcolon}
\global\long\def\dom{\mathop{\mathrm{dom}}}
\global\long\def\argmax{\mathop{\mathrm{arg\,max}}}
\global\long\def\argmin{\mathop{\mathrm{arg\,min}}}
\global\long\def\supp{\mathop{\mathrm{supp}}}

\global\long\def\transp{\top}
\global\long\def\indic{\mathbb{\mathrm{1}}}

\global\long\def\trace{\mathop{\mathrm{Tr}}}
\global\long\def\CVaR{\mathop{\mathrm{CVaR}}}

\global\long\def\DF{\mathrm{D}}
\global\long\def\SD{\text{\ensuremath{\partial}}}

\global\long\def\colplus{\oplus}
\global\long\def\rowplus{\oplus}
\global\long\def\diag{\mathop{\mathrm{diag}}}

\global\long\def\cddot{\mathop{\cdot\cdot}}

\global\long\def\deltaX{\hat{X}}
\global\long\def\deltaY{\hat{Y}}
\global\long\def\deltaZ{\hat{Z}}
\global\long\def\deltaU{\hat{U}}
\global\long\def\deltaB{\hat{B}}
\global\long\def\deltaSigma{\hat{\Sigma}}
\global\long\def\deltaF{\hat{F}}
\global\long\def\deltaPhi{\hat{\Phi}}

\global\long\def\deltaXast{\deltaX^{\ast}}
\global\long\def\deltaYast{\deltaY^{\ast}}
\global\long\def\deltaZast{\deltaZ^{\ast}}
\global\long\def\deltaUast{\deltaU^{\ast}}

\global\long\def\interior{\mathop{\mathrm{int}}}
\global\long\def\sgn{\mathop{\mathrm{sgn}}}
\global\long\def\cvar{\mathbb{ES}}
\global\long\def\var{\mathrm{VaR}}
\global\long\def\as{\text{-a.s.}}

\global\long\def\warrow{\stackrel{\text{w}}{\longrightarrow}}
\global\long\def\wsarrow{\stackrel{\text{w}^{\ast}}{\longrightarrow}}

\global\long\def\cl{\mathop{\mathrm{cl}}}
\global\long\def\diam{\mathop{\mathrm{diam}}}
\global\long\def\ba{\mathrm{ba}}

\global\long\def\rhoMD{\rho^{\mathrm{MD}}}
\global\long\def\rhoMSD{\rho^{\mathrm{MD+}}}
\global\long\def\rhoEpsMSD{\rho_{\epsilon}^{\mathrm{MD+}}}
\global\long\def\rhoEntr{\rho^{\mathrm{Ent}}}

%\KEYWORDS{deterministic inventory theory; infinite linear programming duality; 
%  existence of optimal policies; semi-Markov decision process; cyclic schedule}



%%%%%%%%%%%%%%%%%%%%%%%%%%%%%%%%%%%%%%%%%%%%%%%%%%%%%%%%%%%%%%%%%%%%%%

% Samples of sectioning (and labeling) in IJOC
% NOTE: (1) \section and \subsection do NOT end with a period
%       (2) \subsubsection and lower need end punctuation
%       (3) capitalization is as shown (title style).
%
%\section{Introduction.}\label{intro} %%1.
%\subsection{Duality and the Classical EOQ Problem.}\label{class-EOQ} %% 1.1.
%\subsection{Outline.}\label{outline1} %% 1.2.
%\subsubsection{Cyclic Schedules for the General Deterministic SMDP.}
%  \label{cyclic-schedules} %% 1.2.1
%\section{Problem Description.}\label{problemdescription} %% 2.

% Text of your paper here
\section{Introduction}
Traditionally, the mean of the response is a widely-used performance measure of a stochastic system. However, 
the mean itself is not able to evaluate the possible variability or describe the entire distribution adequately. To provide more thorough profiles of the response distribution, the quantile has become increasingly popular and of great interest in many fields, including insurance, engineering safety, finance, and healthcare \citep{wipplinger2007philippe,morgan1996riskmetrics,cope2009challenges}. In risk management, quantile, also termed as Value-at-risk (VaR), is one of the primary risk measures to quantify and interpret the risk that one system may face. For instance, in the finance industry, the $\alpha$ quantile of a loss distribution represents the lower bound of large losses that the investor can suffer from an activity, where the large losses are defined to be the upper $(1-\alpha)$-tail of distribution with $\alpha$ very close to 1 (like 0.95, 0.99) \citep{hong2009simulating}. 

Optimizing the quantiles of loss functions is a common practice for decision-makers to manage the risk. In this case, searching the best design with the smallest $\alpha$-quantile of $L$ will return the desired decision. More formally, for design choice $x\in \mathcal{X} \subset \mathbb{R}^d$ ($\mathcal{X}$ is the design space assumed compact), we want to minimize the $\alpha$-quantile, $v_{\alpha}(L(x)):=\inf\{y|F_x(y)\geq \alpha \}$, for loss function $L(x)$ (with $F_x(\cdot)$ and $f_x(\cdot)$ defined as the cumulative distribution function and probability density function of $L(x)$):
\begin{equation}
\label{objective}
\text{min}_{x\in \mathcal{X}} \ \  v_{\alpha}(L(x)).
\end{equation}
As large losses are typically of interest, in this work, we consider high quantiles (whose level $\alpha$ is close to 1) of $L(x)$. 

\subsection{Motivation}
The optimization problem \eqref{objective} can be challenging for a few reasons. First, the loss function $L(x)$ usually has no closed-form and is difficult or expensive to observe from the real system. Instead, some simulation engines for $L(x)$ are built, such as the financial model for risk management. Therefore, optimizing $v_{\alpha}(L(x))$ is often conducted via its simulation and a lot of Monte Carlo methods based on the simulation results have been developed (see \citet{hong2014monte} for a review). Second, even with possible simulation models for $L(x)$, $\eqref{objective}$ is still not easy to solve, as $v_{\alpha}(L(x))$ is not directly returned by simulation results but estimated from them. When $\alpha$ is large, it may require a large number of simulations to estimate $v_{\alpha}(L(x))$ precisely. This can be seen from the noise of the quantile estimator. If we denote $n$ as the number of simulation replications at $x$, the noise variance of the empirical quantile estimator is approximately $\alpha(1-\alpha)/[nf_x^2(v_{\alpha}(L(x)) )]$ \citep{bahadur1966note}. For $\alpha$ close to 1, $f_x(v_{\alpha}(L(x)) )$ is typically quite small, especially for heavy-tailed distributions. Therefore, a large number of simulations are required to obtain accurate quantile estimator. The simulation models, however, can be very complicated and time-consuming due to the complex nature of the real system. This restricts the applicable number of simulation runs and makes it almost impossible to obtain results for every considered design with a limited budget. Third, the quantile functions may be non-convex and thus difficult to optimize.  

With these challenges, some optimization algorithms via simulation can be designed to solve \eqref{objective}. The proper algorithm should have at least the following two characteristics. First, it should not require some strict properties from the objective functions, like convexity. Second, it should be efficient and can be used for expensive simulations with a limited budget. In this work, we aim to develop a metamodel-based simulation optimization algorithm which satisfies both characteristics. 

\subsection{Literature Review}
The main idea of metamodel-based simulation optimization approach is to introduce a statistical model to guide the search when optimizing black-box functions. With a limited budget, we can only observe the objective functions at a small number of design inputs, while at the unknown regions, the metamodel serves as an approximation of the true response surface. It provides the information about the entire space and helps decide new points to locate the optimum efficiently. This type of approach has been successfully used in optimizing expensive functions \citep{jones1998efficient,srinivas2009gaussian,regis2007stochastic,muller2017socemo}. It can be classified with respect to the type of metamodel adopted. Some commonly-used metamodels including polynomial regression, radial basis functions, Gaussian process model, artificial neural networks (see \citet{barton2006metamodel} and \cite{jones2001taxonomy} for reviews). Among these methods, the Gaussian process (GP, also termed as kriging) model has become popular as it provides an estimate of the prediction uncertainty, which can be used to construct the selection criterion for further design choices. In this work, we also adopt the GP type metamodel.

Based on the GP model, a few different simulation optimization approaches have been proposed. For deterministic problems (the objective function $v_\alpha(L(x))$ in \eqref{objective} is replaced by some deterministic function $f(x)$), the Efficient Global Optimization (EGO) \citep{jones1998efficient} algorithm with Expected Improvement (EI) criterion is the most widely used for its capability to balance between exploration (searching unexplored region) and exploitation (searching the current promising region). As the function value $f(x)$ can be simulated with no error, no replications at each design input are needed and thus EGO only considers how to select new design points. In parallel with the EGO paradigm, a few different algorithms were developed. GP upper confidence bound (GP-UCB) algorithm provided an alternative to negotiate exploration and exploitation with a tuning hyperparameter making balance between them \citep{srinivas2009gaussian}. Stepwise uncertainty reduction (SUR) approach was to reduce an uncertainty measure with sequentially selected design points \citep{picheny2015multiobjective}. Moreover, some information-based algorithms were developed considering the distributions of the global minimizer (see \cite{shahriari2016taking} for a review). For stochastic simulations, $f$ is simulated with noise and the expected value of $f$ is often considered to be optimized (objective function in \eqref{objective} becomes $\mathbb{E}[f(x,\xi)]$, where $\xi$ represents the randomness). The noise in response needs to be taken care of by the metamodel and optimization algorithm. For response with homoscedastic noises, \citet{huang2006global} proposed Sequential Kriging Optimization (SKO) with the nugget effect GP as metamodels. It introduced an augmented EI to consider the `usefulness' of more replications at one location. For responses with heterogeneous noises, \cite{picheny2013quantile} proposed the Expected Quantile Improvement (EQI), which is an extension of EI, to consider the known noise levels at both the already observed design points and the future candidate. When the noise levels are unknown, recently, some algorithms were proposed including Two Stage Sequential Optimization (TSSO) \citep{quan2013simulation} and extended TSSO (eTSSO) \citep{Pedrielli2018etsso} with the stochastic GP \citep{ankenman2010stochastic,yin2011kriging} as metamodels. They tried to combine the Optimal Computing Budget Allocation (OCBA) \citep{chen2000simulation} technique to decide the number of replications at design point with the EI criterion. The spatial uncertainty of the GP model and the noises of the observations are then reduced iteratively with a global and local search. 

Instead of optimizing the expectation of stochastic functions, in this work, we aim to optimize the quantile functions of loss distributions. Therefore, a metamodel for the quantile function is required. Developing metamodels for quantiles has been extensively studied \citep{koenker2005quantile,dabo2010robust,chen2009metamodels}. Among these models, the quantile regression (QR) \citep{koenker2005quantile} is the primary and most widely used. Recently, the stochastic GP model has been generalized for quantile metamodeling \citep{chen2016efficient}. It shows competitive performance compared with QR model and thus enables us to integrate the GP model into some optimization algorithms for quantile optimization.

\subsection{Illustration \& Contributions}
Using the generalization of \cite{chen2016efficient}, we can extend the eTSSO algorithm for Quantile (eTSSO-Q) optimization. This extended algorithm, however, can still be costly for high quantiles with a limited computing budget. To address this challenge and further improve efficiency, we propose a novel eTSSO-Q Multi-Level (eTSSO-QML) algorithm, which is the main contribution of this work. Different from traditional approaches which directly optimize the quantile function at the objective level, eTSSO-QML starts with optimizing some lower quantiles. Typically, the lower quantiles are cheaper and easier to estimate and their estimations are likely to be less noisy compared with that of a high quantile \citep{bahadur1966note}. 

We next illustrate this idea with an test function from \citet{shim2009non} (the quantile functions are shown in Figure 1).
\begin{figure}[h]
	\centering\includegraphics[width=0.6\linewidth]{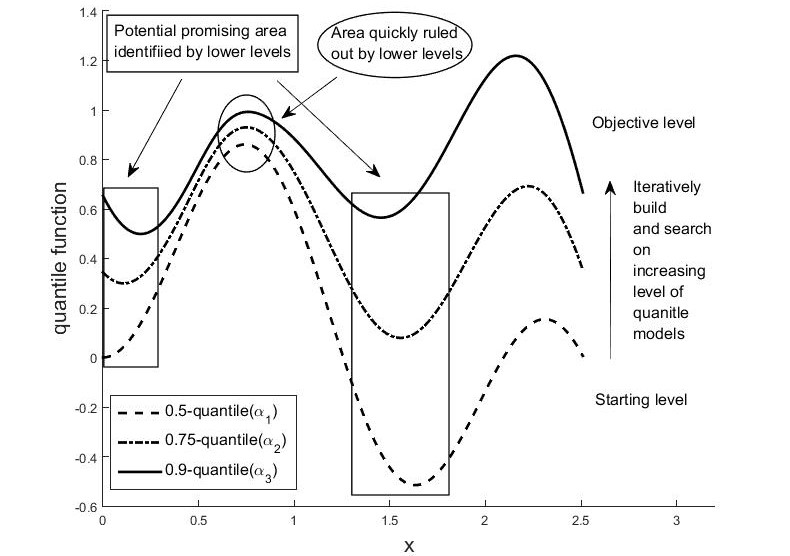}
	\caption{Illustration of eTSSO-QML with  $L(x)=\sin(2.5x)\sin(1.5x)+\mathcal{N}(0,0.01+0.25[1-\sin(2.5x)]^2)$ }
\end{figure}
In the example, we first optimize the more accurate lower levels (0.5, 0.75) to identify promising regions (near 0 and 1.5 in Figure 1). As the algorithm proceeds, more simulation replications are assigned and the quantile estimators at higher levels (0.9) improve. We then increase the level of the quantile being optimized iteratively up to the objective level. At the same time, the search process is guided by the metamodels for these increasing levels with a focus on the possible promising regions identified by the previous levels. In contrast, if we directly optimize the 0.9 quantile, due to the inaccurate quantile estimators obtained with a limited budget, the constructed metamodel can be very unreliable and can mislead the search, resulting in inefficient usage of the budget. 

Optimizing the lower quantile functions can be informative for the objective level for a few reasons. First, the quantile functions at different levels are likely to be correlated since they come from the same loss distribution $L(x)$ \citep{wang2017joint}. In this case, as the quantile level approaches the objective level, the shapes of the surfaces tend to be similar, and thus the solutions found by the previous levels are likely to be promising for the objective level. Second, consider two levels 0.5 and 0.9 in the example above. For $x$ near 0.8, $v_{0.5}(L(x))$ is very large leading us to conclude that the $v_{0.9}(L(x))$ here will be even larger due to the non-decreasing property of the quantile functions (i.e. $v_{\alpha_2}(L(x))\geq v_{\alpha_1}(L(x))$ for $\alpha_2>\alpha_1$). In particular, when we see here that $v_{0.5}(L(0.8))$ is larger than $v_{0.9}(L(1.5))$, it is obvious that $x=0.8$ cannot be optimal for $v_{0.9}$ and hence, there seems no need to allocate further replications to the region near 0.8. In this sense, leveraging the lower levels may help eliminate some bad regions and thus can improve the algorithm efficiency.

More formally, with eTSSO-QML, we consider the problem where the $\alpha_m$-quantile for loss function $L(x)$ is to be minimized leveraging on $m-1$ lower quantiles $0<\alpha_1<\alpha_2<...<\alpha_m<1$. Our main contributions can be summarized as follows:
\begin{enumerate}
	\item We propose a multi-level co-kriging model for the $m$ quantile functions. This model ensures that the predictive curves for different quantiles do not cross and thus the non-decreasing property of quantile functions is maintained. 
	\item With the proposed metamodel, we design the eTSSO-QML algorithm. This algorithm leverages on the multi-level model and starts by searching informative and cheaper lower quantiles to quickly identify promising regions for the objective $\alpha_m$-quantile level.
	\item We prove the convergence of eTSSO-QML and test its empirical performance with several numerical examples.
\end{enumerate}

The rest of this article is organized as follows. Section 2 reviews co-kriging model basics and Section 3 extends it to the multi-level quantile case. Section 4 provides details of eTSSO-QML algorithm and Section 5 states its convergence results. Section 6 provides numerical examples to show the effectiveness of eTSSO-QML. Section 7 summarizes the work and presents some future work. The proofs of all lemmas and theorems are provided in the supplementary material.

\section{A Review of Stochastic Co-Kriging Model Basics}
To jointly model these $m$ quantile functions, we propose to use co-kriging. It was originally developed to model deterministic multi-fidelity problems (where a response can be observed with different fidelities) \citep{kennedy2000predicting} and has recently been extended to stochastic simulation metamodeling for expectations \citep{chen2017stochastic}. In this section, we briefly review some basics of the stochastic co-kriging model.

Here, we first introduce some notations used and the simulation background. To develop a stochastic kriging model, replications of the experiments are required. That is, at each design input $x$, a few simulation runs are required. Throughout this work, we use $L$ to represent the results of the simulations (or equivalently, the random samples of the loss distributions from simulations). For instance, where there are $n$ simulations at $x$, we observe $n$ results: $L(x, \xi_1),...,L(x,\xi_n)$, where $\xi_i$ represents the randomness of the $i$th simulation. With these simulation results, a point estimate for the response of interest (denoted by $\mathcal{Y}(x)$) can be obtained. For instance, when modeling the quantile function, $\mathcal{Y}(x)$ is the sample $\alpha$-quantile: $\mathcal{Y}(x)=L_{\llcorner \alpha n \lrcorner}(x) $, where $L_{\llcorner i \lrcorner}(x)$ is the $i$-th order statistic for the sample $L(x,\xi_1),...,L(x,\xi_n)$. Due to a limited number of simulation runs that can be conducted, $\mathcal{Y}(x)$ is a noisy estimate. Point estimates $\mathcal{Y}(x)$ taken at all design inputs can then be used to develop a predictive model.

The standard stochastic co-kriging model is designed for expectations of a series of stochastic responses. Models at different levels satisfy the following relations:
\[Y_l(x)=Z_l(x)+\epsilon_l(x)=\rho_{l-1}Z_{l-1}(x)+\delta_l(x)+\epsilon_l(x),  \ \ \text{if} \ \  1<l\leq m,  \]
\[Y_l(x)=Z_l(x)+\epsilon_l(x)=\delta_l(x)+\epsilon_l(x), \ \ \text{if} \ \  l=1,\]
where $Y_l$ and $Z_l$ represent the noisy and noise-free responses at level $l$, respectively, and $\delta_l(x)$ ($l=1,...,m$) are $m$ independent second-order stationary GPs \citep{santner2013design}. In each model $\delta_l(x)$, for any finite set of $\{x_1,...,x_t\}$, the GP value $\{\delta(x_1),...,\delta(x_t)\}$ follows a multivariate Gaussian distribution with mean ${{f}}_l(x)^T \beta_l$ and pairwise covariance: $\text{cov}(\delta_l(x_1),\delta_l(x_2))=\sigma^2_l \text{corr}_l (x_1,x_2)$. Here, ${f}_l(x)$ is a $p_l \times 1$ vector of known functions and $\beta_l$ is a vector of model parameters. Without prior knowledge of the mean functions, ${f}_l(x)=1$ is used in this work for illustration. For the correlation function, we adopt the popular Gaussian function: $\text{corr}_l (x_1,x_2)=\exp \left\{ \sum_{j=1}^{d} \frac{{(x_{1,j}-x_{2,j})}^2}{-\theta_{l,j}} \right\} $, 
where $x_{i,j}$ is the $j$th coordinate of $x_i$ and $ \theta_l=(\theta_{l,1},...,\theta_{l,d})$ is the sensitivity parameter determining how large the correlation is in each dimension of $x$. The $m$ random noises, $\epsilon_l$, $l=1,...,m$, follow an $m$-dimensional normal distribution with zero mean. These noises are assumed to be independent of $\delta_l$. It is clear that in this model, $Z_l$ is represented by a scaled $Z_{l-1}$ term, $\rho_{l-1}Z_{l-1}$ plus a difference term. This type of autoregressive structure is first introduced by \citet{kennedy2000predicting} for deterministic multi-fidelity problems.

When the estimates of the responses at some selected design points are obtained, the prediction at any unknown point in $\mathcal{X}$ can be computed based on the co-kriging model. Denote $D$ as the set of design points with $|D|$ representing its cardinality and $\mathcal{Y}_l(x)$ as the point estimate for $Y_l(x)$ for $x\in D$. We assume the design sets for all levels of $Y_l$ are the same and thus for all $x\in D$, the estimates $\mathcal{Y}_l(x)$, $l=1,...m$ are available. This assumption holds in our multi-level quantile case since the point estimators can be obtained by the order statistics of $L(x, \xi_1),...,L_n(x,\xi_n)$ for all desired quantiles. With $\mathcal{Y}^T=(\mathcal{Y}_1(D)^T,...,\mathcal{Y}_m(D)^T)$, where $\mathcal{Y}_i(D)^T:=(\mathcal{Y}_i(x) )_{x\in D}$ is the point estimate vector of the $i$th level for points in $D$, the predictor and its predictive variance of $Z_l(x)$ at any unobserved point $x\in \mathcal{X}$ can be derived as \citep{chen2017stochastic}:

\begin{equation}
\label{eqn:1.0}
\widehat{Z}_l(x):=h_l(x)^T\widehat{\beta}+t_l(x)^T R^{-1} (\mathcal{Y}- H \widehat{\beta}),
\end{equation}

\begin{equation}
\label{eqn:2.0}
\text{var}(\widehat{Z}_l(x)) := \sigma_l^2+ \sum_{j=1}^{l-1}{(P_j^{l-1})}^2\sigma_j^2-t_l(x)^T R^{-1}t_l(x) + \zeta_l(x)^T{(H^TR^{-1}H )}^{-1}\zeta_l(x).
\end{equation}
The notations used in \eqref{eqn:1.0} and \eqref{eqn:2.0} are listed in Table 1. If only one response is considered, $m=1$, we get the stochastic GP model \citep{ankenman2010stochastic}. Furthermore, if the response is observed with no noise, $m=1$ and $R_\epsilon=0$, we get the deterministic GP model.

\begin{table}[h]\small
	\centering
	\caption{Notations list in \eqref{eqn:1.0} and \eqref{eqn:2.0} \label{tab: notation}}
	\begin{tabular}{c|c}
		\hline
		Notation & Definition\\
		\hline
		$P_i^j$ & Products of $\rho_l$. $P_i^j=\prod_{k=i}^{j} \rho_k$ if $j>i$; $P_i^{i-1}=1$  \\
		\hline \multirow{2}{*}{$A_k$}& Correlation of the design points generated by $\delta_k$ whose $(p,q)$-th entry is $A_{k,{pq}}= \text{corr}_k (x_p,x_q)$, \\ 
		& where $x_p$ and $x_q$ are the $p$th and the $q$th design point in $D$, respectively.\\ \hline
		$A_j(D,x)$ & Correlation between $x$ and the design points generated by $\delta_j$ \\ \hline
		$\zeta_l(x)$ & $h_l(x)^T-t_l(x)^T R^{-1}H$ \\ \hline $h_l(x)^T$ &$(P_1^{l-1}{\bf{f}}_1(x)^T, P_2^{l-1}{\bf{f}}_2(x)^T,..., \ \ P_{l-1}^{l-1}{\bf{f}}_{l-1}(x)^T, {\bf{f}}_{l}(x)^T, {\bf{0}}^T_{p_{l+1}+...+p_m})$ \\ \hline $t_l(x)^T$ &$(t_{l,1}(x)^T,...,t_{l,m}(x)^T)$, with
		$ t_{l,s}(x)= \sum_{j=1}^{q}\sigma_j^2 P_j^{s-1}P_j^{l-1}A_j(D,x), \ \ \text{where} \ \ q=\text{min}\{s,l\}. $\\ \hline \multirow{2}{*}{$R_z$} & The covariance matrix of the spatial uncertainty \\ & A symmetric matrix with $m\times m$ blocks:
		$R_z^{(k,s)}=\sum_{j=1}^{q}\sigma_j^2 P_j^{k-1}P_j^{s-1}A_j, \ \ \text{where} \ \ q=\text{min}\{k,s\}. $ \\  \hline \multirow{2}{*}{$R_\epsilon$} & The covariance matrix of the noises \\ & A symmetric matrix with $m\times m$ blocks:
		$R_\epsilon^{(k,s)}=  \text{diag}(\text{cov}(\epsilon_k(x_1),\epsilon_s(x_1) )   ,..., \text{cov}(\epsilon_k(x_{|D|}),\epsilon_s(x_{|D|}) ))  .$ \\ \hline $R$ & $R_z+R_\epsilon$ \\ \hline $H$ & A matrix with $m\times m$ blocks: $H^{(k,s)}=P_s^{k-1} {\bf{f}}_s(D)^T,$ if $k \geq s$; $H^{(k,s)}={\bf 0}_{|D| \times p_s },$ if $k < s. $ \\ \hline
		$\widehat{\beta}$ & Best linear unbiased estimator for $\beta$: ${(H^TR^{-1}H)}^{-1}H^TR^{-1}\mathcal{Y}$  \\
		\hline
	\end{tabular}
\end{table}

The above results assume known hyperparameters $\rho_l, {\theta}_k, \sigma_k^2$, $l=1,...,m-1,k=1,...,m$, and covariance matrix for noise, $R_\epsilon$. When building the model in practice, these are typically unknown and should be estimated. Depending on how they are estimated, we separate these hyperparameters into two categories: model inputs and model parameters. The model inputs include the point estimates vector $\mathcal{Y}$ and the estimators for the associated noise covariance matrix $R_\epsilon$. These estimators are directly drawn from the initial simulation results and serve as the inputs to the co-kriging model. For the standard stochastic co-kriging model for the mean performance measures, the inputs are the sample means and sample covariance for the mean estimates. The remaining hyperparameters (${\bf \rho}=(\rho_1,...,\rho_{m-1}),\theta=(\theta_1,...,\theta_m), \sigma^2=(\sigma_1^2,...,\sigma^2_m)$) are referred to as model parameters and can be estimated by maximizing the likelihood function for point estimate vectors (see Appendix A for the likelihood function and some detailed discussion). After this, the predictor (\ref{eqn:1.0}) and predictive variance (\ref{eqn:2.0}) can be obtained by plugging in the estimated parameters.

\section{Stochastic Co-kriging Model for Quantiles}
When applied in quantile predictions, the predictive model structures remain the same as \eqref{eqn:1.0} and \eqref{eqn:2.0}. However, several important adaptions are required. First, we need to find proper approaches to estimate the model inputs, $\mathcal{Y}$ and $R_\epsilon$, which are the point estimate and noise covariance matrix for quantiles instead of expectations in traditional co-kriging model. Section 3.1 introduces the estimation of these inputs and derives some of their properties. Furthermore, due to the non-decreasing property of quantiles, the predictive curves for different levels of quantiles should not cross (which is a criterion not considered in traditional co-kriging models). In Section 3.2, we propose a penalized maximum likelihood estimation (PMLE) approach to ensure non-crossing of our estimates.

\subsection{Estimation of Model Inputs} 

In practice, $\mathcal{Y}$ and $R_\epsilon$ are calculated from simulation results and then plugged into \eqref{eqn:1.0} and \eqref{eqn:2.0}. Specifically, given original simulation results $L(x, \xi_1),...,L_n(x,\xi_n)$ at $x$, we can easily obtain the point estimates for $v_{\alpha_j}(L(x))$ and $v_{\alpha_k}(L(x))$:
\[\mathcal{Y}_j(x)=L_{\llcorner \alpha_j n \lrcorner}(x), \ \ \mathcal{Y}_k(x)=L_{\llcorner \alpha_k n \lrcorner}(x).\]
Following the recommendations of \citet{chen2016efficient} who tested different approaches to estimate the noise variance of the quantile estimates applied in the GP model, including batching \citep{seila1982batching}, sectioning \citep{asmussen2007stochastic}, sectioning-batching \citep{nakayama2014confidence} and jackknifing \citep{nance2002perspectives}, here we use the sectioning method to estimate $\text{var}(\mathcal{Y}_j(x))$ and $\text{var}(\mathcal{Y}_k(x))$. Furthermore, as the noise of these two estimates are correlated since they are drawn from similar simulation results, in this section, we also propose a sectioning method to estimate this noise covariance and derive the asymptotic properties of this estimator.

With the sectioning method, the $n$ simulation runs are first divided into $n_b$ batches with $n_c$ runs in each batch ($n=n_b\cdot n_c$). Then the covariance of ${ {\mathcal{Y}}_{j,l}(x)}$ and ${ {\mathcal{Y}}_{k,l}(x)}$ is estimated based on the quantile estimators with all simulation runs at $x$, $\mathcal{Y}_j(x)$ and $\mathcal{Y}_k(x)$, and the estimators within each batch, ${ {\mathcal{Y}}_{j,l}(x)}$ and ${ {\mathcal{Y}}_{k,l}(x)}$, where ${ {\mathcal{Y}}_{j,l}(x)}$ and ${ {\mathcal{Y}}_{k,l}(x)}$ are the sample $\alpha_j$ and $\alpha_k$ quantiles of the $l$th batch: $L(x,\xi_{1}^l),...,L(x,\xi_{ n_c}^l)$, $l=1,...,n_b$, respectively.
\begin{equation}
\label{noiseva}
\widehat{\text{var}}(\epsilon_i(x))=\widehat{\text{var}}(\mathcal{Y}_i(x))=\frac{1}{n_b(n_b-1)}\sum_{l=1}^{n_b}{({ {\mathcal{Y}}_{i,l}(x)}-{ {\mathcal{Y}}_{i}(x)})}^2, 
\end{equation}
\begin{equation}
\label{noiseco}
\widehat{\text{cov}}(\epsilon_j(x),\epsilon_k(x)) =\widehat{\text{cov}}({ {\mathcal{Y}}_{j}(x)},{ {\mathcal{Y}}_{k}(x)}) =\frac{1}{n_b(n_b-1)}\sum_{l=1}^{n_b}({ {\mathcal{Y}}_{j,l}(x)}-{ {\mathcal{Y}}_{j}(x)})({ {\mathcal{Y}}_{k,l}(x)}-{ {\mathcal{Y}}_{k}(x)}),
\end{equation}

\begin{assumption}
	For all $x\in \mathcal{X}$, $L(x)$ has continuous distribution $F_x$ with density function $f_x$, and finite mean and variance. The function $f_x$ has bounded first order derivatives in the neighborhood of $v_\alpha(L(x))$ with $f(v_\alpha(L(x)))>0$, where $v_\alpha(L(x))$ is the true $\alpha$-quantile.
\end{assumption}

Under Assumption 3.1, \cite{chen2016efficient} has shown that $\widehat{\text{var}}(\mathcal{Y}_k(x))$ is asymptotically unbiased with mean squared error (MSE) of order $o(n^{-2})$ as $n_b,n_c \rightarrow \infty$. Following a similar approach, we can also prove the asymptotic properties of the proposed noise covariance estimator \eqref{noiseco}.
\begin{thm}
	Under Assumption 3.1, when $n_b,n_c\rightarrow \infty$, $\widehat{\text{cov}}({\mathcal{Y}}_{j}(x),{\mathcal{Y}_{k}}(x))$ is consistent and asymptotically unbiased, and the MSE of $\widehat{\text{cov}}({\widehat {\mathcal{Y}}_{1}}(x),{\widehat {\mathcal{Y}}_{2}}(x))$ is $o(n^{-2})$.
\end{thm}

\subsection{A PMLE Approach to Avoid the Crossing Problem}

Traditional QR models quantile functions at different levels separately, which can result in possible crossing between different quantile predictive curves. This, for example, will cause the predictive value of the 0.95 quantile at some points to be larger than that of the 0.99 quantile. This crossing phenomenon is a widely acknowledged problem in quantile modeling and can lead to an invalid distribution of the response and problematic inferences \citep{koenker1984note,cole1988fitting,he1997quantile}. For our quantile co-kriging model, preventing crossing to ensure monotonicity not only improves inferences but more importantly ensures that the multi-level search in the optimization algorithm is valid and efficient. Imagine if the crossing happens between two quantile models, the non-promising region identified by the lower quantile model can be misleading, since the higher quantile can be smaller than the lower one, and hence, can have promising (and even optimal) values in those non-promising regions. Therefore, it is vital for model validity and optimization efficiency to ensure non-crossing in the models before the optimization process. Although in the co-kriging model, multiple quantiles are modeled jointly, non-crossing is not guaranteed. Note that in the traditional application of the co-kriging model where deterministic or mean responses have typically been modeled, the crossing of the models is not a problem.

In this section, we propose a new penalized version of the stochastic co-kriging model to prevent crossing for the quantile models. In our multi-level quantile problem, there is no crossing when:
$$\widehat{Z}_{l+1}(x)-\widehat{Z}_l(x) \geq 0, \ \ \forall x\in \mathcal{X}, l=1,...,m-1.  $$
In other words, the difference between the predictive curves for any two successive quantiles should be non-negative across the design space. We first propose a penalized GP model that can ensure non-negative predictions for a single deterministic response (which can be considered as the difference between two quantiles and thus is non-negative everywhere in the design space) and then apply it to our multi-level model.

Consider first a deterministic GP model for a non-negative function $W$ (to distinguish with the model in the previous section, we use $W$ here to represent this response and $\mathcal{W}$ as the observations for it):
\begin{equation}
\label{Wmodel}
W(x)={f}(x)^T \beta + M(x;\theta) ,
\end{equation}
where ${f}(x)$ is a $p \times 1$ known function, $\beta$ is a $p \times 1 $ vector of model parameters and $M(x;\theta)$ is assumed to be a zero-mean second-order-stationary GP controlled by hyperparameters $\theta$. Given that the true function $W$ is non-negative and that the observations have no noise, the observation vector we get, $\mathcal{W}$, is non-negative. The standard GP model for a deterministic function is actually an interpolation of the observations $\mathcal{W}={(\mathcal{W}(x_1),...,\mathcal{W}(x_t))}^T$, and the shape of the predictive curve changes with the hyperparameters $\theta$. Therefore, when estimating $\theta$, we must make sure that the resulting curve should not intersect with the surface $W=0$. In other words, those values of $\theta$ that cause the intersection should be eliminated. This intuition can naturally translate into the following penalization method. Instead of optimizing the ordinary loglikelihood of $\mathcal{W}$, we propose to minimize the following penalized likelihood function to get the PMLE estimator for $\theta$:
\[
Q(\mathcal{W},\theta):=-l(\mathcal{W},\theta)+P(\mathcal{W},\theta)   =\frac{1}{2} \ln(|{(R')}|)+ \frac{1}{2} {(\mathcal{W}- F \widehat{\beta'})}^T{(R')}^{-1}(\mathcal{W}- F \widehat{\beta'}) + \lambda \cdot \kappa(\mathcal{W},\theta),
\]
where $l(\mathcal{W},\theta)$ is the ordinary loglikelihood function, $P(\mathcal{W},\theta)=\lambda \cdot \kappa(\mathcal{W},\theta)$ is the penalty term, $\lambda$ is a non-negative penalty coefficient, $F={({\bf f}(x_1), ..., {\bf f}(x_n))}^T$, $R'$ is the covariance matrix for $\mathcal{W}$, $\widehat{\beta'}={(F^TR^{-1}F)}^{-1}F^TR^{-1}\mathcal{W}$ and
$$\kappa(\mathcal{W}, \theta)=\left\{\begin{array} {cc} |\text{min}_{x\in \mathcal{X}}(\widehat{W}(x))|, & \text{if}\ \ \text{min}_{x\in \mathcal{X}}(\widehat{W}(x))<0 \\ 0 & \text{if}\ \ \text{min}_{x\in \mathcal{X}}(\widehat{W}(x)) \geq 0 \end{array} \right. .$$ 
With this penalty term, the parameters that cause the predictive curve to go below the $W=0$ plane will be penalized. Theorem 2 demonstrates the consistency of the parameter estimated with this approach. It is established based on the asymptotic property of the MLE for the ordinary GP model (denoted as $\widehat{\theta}_o$). Specifically, under certain regularity conditions, $\sqrt{n}(\widehat{\theta}_o-\theta_0)\rightarrow\mathcal{N}(0, \mathbb{I}^{-1}(\theta_0)) $ as $n\rightarrow \infty$ in distribution, where $n$ is the number of design points and $\mathbb{I}$ is the Fisher information matrix \citep{mardia1984maximum}.
\begin{thm}
	Denote $\theta_0$ as the true value of $\theta$ for model \eqref{Wmodel} and $n$ as the number of design points. There exists a local minimizer $\widehat{\theta}_n$ of $Q(\mathcal{W},\theta)$ such that $||\widehat{\theta}_n-\theta_0||=\mathcal{O}_p(1/\sqrt{n})$.
\end{thm}

This PMLE approach involves an optimization problem over the predictive surface $\widehat{W}(x)$. We highlight that this optimization is much easier compared to optimizing the true unknown response surface since the predictive response function is much cheaper with explicit form. When applied in our case, where the function $W$ (the difference between two quantile functions) is stochastic, this method can also return non-negative predictions by preventing crossing between the predictive curves and the plane $W=0$. With a slight modification of the penalty function $\kappa$, this method can be easily applied in our multi-level quantile problem:
$$\kappa(\varphi)=\left\{\begin{array} {cc} |\varphi|, & \text{if}\ \ \varphi<0 \\ 0, & \text{if}\ \ \varphi \geq 0 \end{array} ,  \right. $$
where $\varphi=\text{min}_{x\in\mathcal{X},l\in \{1,...,m-1\}}(\widehat{Z}_{l+1}(x)-\widehat{Z}_{l}(x))$. It is easy to see that the parameters will be penalized once crossing happens between any two successive predictive curves among the $m$ quantile models.

With the approaches proposed here, we can build a co-kriging model for multi-level quantiles that does not cross. As mentioned before, there exist quite a few different approaches to do this more rigorously. For the GP based model, some other more complicated and refined methods have also been proposed to ensure positive response prediction \citep{szidarovszky1987kriging,dowd1982lognormal}. Compared with those methods, the PMLE approach keeps the nice auto-regressive structure and is convenient to use and integrate into the multi-level algorithm. From a more pragmatic viewpoint, as the metamodel here is mainly used as an aid to the optimization process, we do not consider more sophisticated techniques and just apply the PMLE approach.

\section{Multilevel Quantile Optimization (eTSSO-QML) Algorithm}
This section presents the eTSSO-QML algorithm, which optimizes the $\alpha_m$ quantile with a multi-level model built from the $\alpha_1,..,\alpha_m$ quantiles. As previously noted, the optimization process is guided by the proposed stochastic quantile co-kriging model. It starts with searching the lower quantiles and then searches on the promising regions for higher quantiles identified by the lower ones. The algorithm is fully sequential where the overall computing budget is iteratively allocated. Within each iteration, we apply the two-stage framework from the eTSSO algorithm to provide a ``division of labor" \citep{Pedrielli2018etsso}. In the first stage (Searching Stage), we adopt the EI criterion to select a new design input with the highest probability of achieving a better result than the current best. The second stage (Allocation Stage) focuses on distributing additional simulation replications to the existing design points. This is to improve the model fit and increase our confidence in the estimators to correctly identify the optimum. The distribution of budget used in these two stages is allowed to change with iteration. At the beginning of the algorithm, as little is known about the response, more budget is used to search the design space to identify the promising regions; and towards the end, more budget will be saved for the allocation stage, since the emphasis then becomes refining the point estimates at already selected designs when we are in proximity to the promising regions. These two stages will be discussed in detail in Section 4.3 and Section 4.4 after an overview of the algorithm is given in Section 4.1 and an introduction of the algorithm parameters is given in Section 4.2. Our algorithm is based on the eTSSO procedure, and we refer interested readers to \cite{Pedrielli2018etsso} for full details of the algorithm.

\subsection{Algorithm Overview}

Before describing the algorithm, we list key parameters in Table 2.
\begin{table}[h]\small
	\centering
	\caption{Algorithm parameters list \label{tab: first}}
	\begin{tabular}{cc}
		\hline
		Parameter & Definition\\
		\hline
		$T$ & Total number of replications (computing budget) at the beginning \\ $D_0$& Initial design set\\ 
		$\alpha_l,l=1,...,m$ & Quantiles used for modeling \\$r_0$ &Minimum number of replications for a newly selected design input \\ $C_0$ & The maximum noise variance of a quantile estimate that can be tolerated \\$k$ &Current iteration\\ $D_k $& Design set at iteration $k$\\$h(k)$ & Current level of quantile guiding the search \\ $B_k$& Number of available replications in iteration $k$\\  $\pi_k$ & The set of the quantile levels building the co-kriging model in iteration $k$ \\
		$E_l$ & Set of inputs whose estimates at the $l$th level have acceptable accuracy \\ $\mathcal{Y}_l$ & Observations for $l$th level \\ $A$ & Remaining number of replications (Algorithm terminates when $A=0$)	\\	 $\widehat{x}_k$ & The best input for the $\alpha_m$th quantile found by iteration $k$ \\
		\hline
	\end{tabular}
\end{table}
The first five parameters are user-defined to start the algorithm. The total number of replications, $T$, is typically determined by the computing budget and for $D_0$, if no prior knowledge or preference is available, users can apply non-informative design strategies such as the uniform and Latin Hypercube sampling strategies. The values of $\alpha_1,...,\alpha_m$ should also be specified in advance. As noted above, to optimize a high quantile $\alpha_m$, we start with the base level $\alpha_1$. This level should not be too high and we suggest $\alpha_1\in[0.5,0.6]$ based on our experience. For the remaining level, we consider fixed $(m-2)$ evenly distributed inter-levels between $\alpha_1$ and $\alpha_m$. The number of levels, $m$, can be selected depending on the budget. A larger $m$ can slow the approach to the objective level and increase the co-kriging model complexity. However, as the difference between any two successive levels becomes smaller, the correlation between them increases, and thus the promising regions identified by lower levels become more accurate. In contrast, a smaller $m$ can reduce the computational burden but may weaken the correlation among the levels adopted. For any newly selected design point, we first assign $r_0$ replications to it. This $r_0$ can be chosen to ensure that the point estimates for the base level have reasonable accuracy. The parameter $C_0$ is used to examine the accuracy of a quantile estimator and only the estimator whose variance is smaller than $C_0$ is accepted. These two parameters can be chosen through a cross-validation test over $D_0$.
To achieve this, we can start with a small number of $r_0$, and iteratively increase its value until the model for $\alpha_1$ built with $D_0$ and $r_0$ replications at each input passes the cross-validation test. After that, $C_0$ can be selected as the maximum of the $\alpha_1$ quantile estimators from the points in $D_0$. The other algorithm parameters are updated with each iteration, and these will be described in detail in Section 4.2. 

The eTSSO-QML algorithm is a iterative algorithm, iterating between the Searching Stage and the Allocation Stage until the computing budget runs out. We illustrate the overall idea of the algorithm with the example in Figure 1. At the start of the algorithm, a small budget is first applied. With a small number of replications, the point estimates of the target quantile (0.9) can be inaccurate with high uncertainty (noisy). At this stage, a more reliable lower quantile $\alpha_1$ model (0.5) is first built and used to guide the initial search. In other words, we optimize the first level as a start to identify possible promising regions (like the regions near 0 and 1.5). As the algorithm proceeds, more budget is allocated and the accuracy of the higher quantile estimators improves. The algorithm will then stepwise increase the level of the quantile metamodels developed, and use the current highest level to guide the search. As a result, the algorithm gradually optimizes higher and higher quantile levels with a focus on the promising regions identified by previous levels, to finally optimize the quantile at the target level $\alpha_m$. The eTSSO-QML algorithm is described in Algorithm 1. In Section 4.3 and Section 4.4, we describe in further detail about the Searching Stage and the Allocation Stage.
\begin{algorithm}[h]\small
	\hspace*{0.02in}{\bf Input:} $T$, $D_0$, $r_0$, $\{\alpha_1,...,\alpha_m\}$
	\begin{algorithmic}[1]
		\caption{eTSSO-QML algorithm}
		\State {\bf Initialization}:
		\State $k=0$; $B_1=r_0$; $A =T-|D_0|\times r_0$
		\State Get the simulation results for $D_0$ with $r_0$ replications for each input
		\State Estimate $\mathcal{Y}_1$ for $D_0$, let $h(0)\leftarrow 1$
		\State Fit single stochastic GP model $\widehat{Z}_1(x)$ for $\mathcal{Y}_1$, resulting the predictive random variable $\tilde{Z}_1(x)$
		\State Let $z_0^*=\min_{x\in D_0}(\widehat{Z}_1(x))$
		\While{$A>0$}
		\State {\bf Searching Stage}:
		\State $x_{k+1}=\arg\max_{x\in \mathcal{X}\setminus {D_k}} \mathbb{E}[\max\{z_{\alpha_{h(k)}}^*-\tilde{Z}_{\alpha_{h(k)}}(x),0\}]$
		\State Run $r_0$ replications at $x_{k+1}$ to obtain the quantile estimators and set $D_{k+1}\leftarrow D_k\cup \{x_{k+1}\}$ 
		\State {\bf Allocation Stage}:
		\State Update $B_k$ 
		\State Allocate budget to ensure that each design point has at least $r_k$ replications.
		\State Use OCBA to allocate the remaining replications to selected inputs and run new simulations correspondingly
		\State {\bf Modeling Update}
		\State Set $E_l=\emptyset, l=1,...,m$.
		\State For each selected design $x_i$, decide $l^*=\arg\max \{l| l\in \{1,..m\},\widehat{\text{var}}(\mathcal{Y}_l(x_i))\leq C_0\} $. Set $E_l=E_l\cup \{x_i\}$
		\State Set $h(k+1)$ as the largest value in $\{1,2,...,m\} $ such that $E_{h(k+1)}\neq \varnothing$ 
		\State For each $j \in [1,h(k+1)]$, if $\{l|j<l\leq h(k+1), E_l=E_j \} = \emptyset$, set $\alpha_j \in \pi_{k+1}$
		\State Fit the co-kriging model for the quantile levels in $\pi_{k+1}$ resulting the predictive random variable $\tilde{Z}_{h(k+1)}(x)$ and Let $z_{k+1}^*=\min_{x\in D_h(k+1)}(\widehat{Z}_{h(k+1)})$.
		\State Report $\widehat{x}_{k+1}=\arg\min_{x\in D_k}(\mathcal{Y}_m(x))$ as the current best point and $ \mathcal{Y}_m(\widehat{x}_{k+1})$ as the best solution at the objective level found by iteration $k+1$.
		\State $A \leftarrow A-B_k$, $k \leftarrow k+1$
		\EndWhile
		\State \Return $\widehat{x}=\arg\min_{x\in D_k}(\mathcal{Y}_m(x))$
	\end{algorithmic}
\end{algorithm}

\subsection{Modeling Update in Each Iteration}

In each iteration, the model \eqref{eqn:1.0} and \eqref{eqn:2.0} and the algorithm parameters are updated. In iteration $k$, the search is guided by the $h(k)$-th level, which is the level to be optimized. The value of $h(k)$ gradually increases from 1 to $m$. We choose its value as follows. For each $x$ in $D_k$, we find the largest value $l$ such that $\mathcal{Y}_{l}(x)$ drawn from simulation results has noise variance smaller than $C_0$. After that, we set $x\in E_j$, $1\leq j\leq l$. As a result, $E_j$ consists of design points that have acceptable accuracy at level $j$. After that, we select $h(k)$ as the largest value in $\{1,2,...,m\} $ such that $E_{h(k)}\neq \varnothing$. In this case, we choose $h(k)$ as the current highest level and then we can build a multi-level model for $\alpha_1,...,\alpha_{h(k)}$ in iteration $k$. However, the increasing value of the current highest level $\alpha_{h(k)}$ naturally increases the model complexity and so we would like to select some but not all from the $\alpha_1,...,\alpha_{h(k)}$ quantiles to build the co-kriging model. 

In fact, as the algorithm proceeds, some inter-level quantiles become redundant. Consider when the objective level is 0.95 quantile and we have inter-levels 0.9 and 0.8 quantiles which have similar design sets with acceptable accuracy, we may remove the 0.8 quantile as the 0.9 quantile is closer to our objective. This removal can be partly interpreted by the auto-regressive structure of the co-kriging model: the 0.95 quantile model depends on the previous levels only through the nearest level, the 0.9 quantile. Therefore, the 0.8 quantile function has no contribution if the 0.9 quantile is reasonably good and thus it can be removed. In this sense, we select a subset $\pi_k$ from $\{\alpha_1,...,\alpha_{h(k)} \}$ by removing some redundant inter-levels and building a co-kriging model for levels in $\pi_k$. Specifically, for $\alpha_j$, $1\leq j \leq h(k)$, only when $\{l|j<l\leq h(k), E_l=E_j \} = \emptyset$, we set $\alpha_j \in \pi_k$. Selecting $\pi_k$ in this way also ensures that asymptotically, we only select $\alpha_m$ to build a single stochastic GP model. This is intuitive since when the number of iterations assigned to each design point tends to infinity, the $\alpha_m$-quantile estimators become accurate and there is no need to leverage on the information from lower levels (see Section 5 for detailed discussion). 

The budget $B_k$ changes with $k$ and a specific choice will be introduced in detail in Section 4.4. The values of $ \mathcal{Y}_l$ and $A$ can be easily updated after the two stages finish. At the end of each iteration, we choose the observed lowest quantile value at the objective level $\alpha_m$ as the optimum found by iteration $k$.

\subsection{Searching Stage}
In the Searching stage, we select the next design input based on the following EI criterion \citep{jones1998efficient}:
\[  
\begin{split}
x_{k+1} & = \arg\max_{x\in \mathcal{X}} \mathbb{E}[\max\{z_{\alpha_{h(k)}}^*-\tilde{Z}_{\alpha_{h(k)}}(x),0\}]\\ &= \arg\max_{x\in \mathcal{X}} \left\{ \widehat{s}_{h(k)}(x)\phi(\frac{z_{\alpha_{h(k)}}^*-\widehat{Z}_{\alpha_{h(k)}}(x)}{\widehat{s}_{h(k)}(x)})+(z_{\alpha_{h(k)}}^*-\widehat{Z}_{\alpha_{h(k)}}(x) )  \Phi(\frac{z_{\alpha_{h(k)}}^*-\widehat{Z}_{\alpha_{h(k)}}(x)}{\widehat{s}_{h(k)}(x)}) \right\} ,
\end{split}
\]
where $z_{\alpha_{h(k)}}^*$ is the lowest value of the predictive responses $\widehat{Z}_{\alpha_{h(k)}}$ at $D_{k}$, $\phi$ and $\Phi$ represent the probability and cumulative distribution functions of the standard normal distribution, respectively, and $\tilde{Z}_{\alpha_{h(k)}}(x)$ is a Gaussian random variable with distribution $ \mathcal{N}(\widehat{Z}_{\alpha_{h(k)}}(x),\widehat{s}_{h(k)}^2(x)  )$, where
\begin{equation}
\label{eqn:3.0}
\widehat{s}_l^2(x):=\sigma_l^2+ \sum_{j=1}^{l-1}{(P_j^{l-1})}^2\sigma_j^2-t_l(x)^T R_z^{-1}t_l(x) + \zeta_l(x)^T{(H^TR_z^{-1}H )}^{-1}\zeta_l(x).
\end{equation}
Compared to \eqref{eqn:2.0}, \eqref{eqn:3.0} uses $R_z$ instead of $R$ and thus it only considers the response covariance generated by the spatial process and ignores the noise variance. The rationale behind this is that the Searching Stage is by design to choose new points to reduce the spatial uncertainty (the noise is taken care of by the Allocation Stage). Moreover, \eqref{eqn:3.0} ensures that the EI values at all selected design points are zero so that they will not be reselected in the future iterations. With this criterion, we use the current highest level, $h(k)$, to guide our search. The EI criterion selects the point which has the largest expected improvement with respect to the current best. Typically, the points with small response predictions (from current promising regions) or large predictive variances (from less explored regions) have large EI values. By selecting these points, EI balances between exploitation and exploration. At this new design point, we run $r_0$ simulations, add it to the design set and then update the point estimate vectors accordingly.

\subsection{Allocation Stage} 

In this stage, we adopt the OCBA technique to allocate computing resources to the selected design points. The original OCBA technique, however, is designed for ranking and selection problems with finitely many alternatives. When the number of possible alternatives is infinite, to guarantee the convergence of the algorithm, we make the following assumption on the allocation rule.
\begin{assumption}
	Suppose there exists a sequence $\{r_1,...,r_k \}$ such that $r_{k+1}\geq r_k$, $r_k \rightarrow \infty$ as $k\rightarrow \infty$ and that $\sum_{k=1}^{\infty}\frac{k}{r_k}<\infty$. Denote $N_k(x)$ as the cumulative number of replications assigned to the selected design point $x$ by iteration $k$. It follows that $\min_{x\in D_k}N_k(x)\geq r_k$ for all $k$.
\end{assumption}
This assumption has been used for problems with discrete but infinite alternatives (\cite{hong2006discrete}, \cite{andradottir2006simulation}). Although we consider optimization problems within a continuous domain, this assumption is important to ensure convergence (see the detailed discussion in Section 5). To fulfill this assumption, in the Allocation Stage of eTSSO-QML, we first spare some budget to ensure that for each selected design point, there are at least $r_k$ replications assigned to it by iteration $k$ (including the newly selected input in iteration $k$). After this initial stage, we adopt the OCBA technique to allocate the remaining replications. 

As noted before, the budget of Allocation Stage increases with $k$ to refine the point estimates at selected design inputs. This is intuitive since at the beginning, more budget can be used to search the design space and when $k$ gets larger, we are more likely to be in proximity of the promising region. At this point, we can reduce the number of newly selected designs and assign more budget to the already sampled points in the promising region to refine our estimate of the optimum. We, therefore, let $B_1=r_0$ and $B_k$ increases with iteration and update its value as follows when $k>1$:
$$B_k=\max\{\sum_{i=1}^{|D_0|+k} \max\{0, r_k-N_k(x_i)\}, \lfloor B_{k-1}(1+\frac{\max\limits_{x_i\in D_k } \widehat{\text{var}}({\mathcal{Y}_{h(k)}}(x_i)) }{\max\limits_{x_i\in D_k } \widehat{\text{var}}({\mathcal{Y}_{h(k)}}(x_i))+ \widehat{s}^2_{h(k)}(x_{k+1})} \rfloor\}, $$
where $\widehat{\text{var}}({\mathcal{Y}_{h(k)}}(x_i))$ is the sample noise variance of the point estimate at $x_i$ (estimated by \eqref{noiseva}). This adaptive scheme $B_k$ is first adopted in the eTSSO algorithm. Its increase is controlled by the relationship between the point estimator noise, measured by $\max_{x_i\in D_k } \widehat{\text{var}}({\mathcal{Y}_{h(k)}}(x_i))$, and the spatial uncertainty of the GP model, measured by the predictive variance $\widehat{s}^2_{h(k)}(x_{k+1})$ (see equation \eqref{eqn:3.0}). At the beginning, when the spatial uncertainty is very large, $B_k$ has a slow growth to save more budget for new design point selection. When the spatial uncertainty gets smaller, i.e., the design space has been better explored, $B_k$ will then experience a faster growth, focusing more on the already selected points in the promising regions. The advantage of this scheme is that it does not require user-defined budgets for each iteration. Furthermore, it can improve the identification of the optimum and lead to efficient use of the computing budget.

After updating $B_k$ and checking that each design point has at least $r_k$ replications, we can allocate remaining replications to the selected design inputs with the OCBA technique. Denote $x_b$ as the current best design point in $D_k$ with respect to the current highest level: $x_b \in \arg \min_{x\in D_k} \mathcal{Y}_{h(k)}(x)$. OCBA decides the number of new replications $n_i$ assigned to each input $x_i\in D_k$ as follows:
$$\frac{n_i}{n_j}=\frac{  \sqrt{\widehat{\text{var}}({\mathcal{Y}_{h(k)}}(x_i))}/{\lambda_{b,i}} }{\sqrt{\widehat{\text{var}}({\mathcal{Y}_{h(k)}}(x_j))}/{\lambda_{b,j}} } , \ \ (i,j \neq b); \ \ n_b= \sqrt{\widehat{\text{var}}({\mathcal{Y}_{h(k)}}(x_b))} \sqrt{\sum\limits_{i\neq b}\frac{n_i}{\widehat{\text{var}}({\mathcal{Y}_{h(k)}}(x_i))}},$$
where $\lambda_{b,i}:=\mathcal{Y}_{h(k)}(x_i)-\mathcal{Y}_{h(k)}(x_b)$. The OCBA technique actually prefers allocating additional replications to points with low response values and large noises, which is intended to refine the point estimates at promising regions and those with large noise variances. Next, we can run additional simulations for the existing inputs and update the point estimate vector, $\mathcal{Y}$, and current highest level, $h(k)$, accordingly. Finally, a new co-kriging metamodel can be built and then the algorithm goes to the next iteration.

As the overall algorithm proceeds, the search focuses more and more on promising regions of higher quantiles, and as a result, less budget is spent in the non-promising regions. This can be seen through our two-stage procedures. In the Searching Stage, as the lower quantiles are easy to estimate, it generally takes a small budget to search the lower quantile. Through this search, the algorithm has a focus on the promising regions for lower levels. If these regions are also promising for $\alpha_m$, it would have already been sampled. In the Allocation Stage, the budget used for increasing the precision of point estimate at lower quantiles essentially also improves the estimates for higher quantiles. Therefore, when optimizing lower quantiles, more budget is spent in its promising regions and the estimates (for all levels of quantiles) are improved. In other words, the algorithm digs into the promising regions identified by the lower quantiles when estimating $\alpha_m$. We highlight that this improvement is due to the specificity of our problem and does not apply to general multi-fidelity problems. In those problems, the experiments for lower fidelity and higher fidelity are different. Running simulations for one level does not necessarily improve the point estimates for another level. 

Before we close this section, we briefly analyze how to make the current highest level, $h(k)$, approach the objective level, $m$, i.e., to ensure that the objective level is optimized. This can be easily achieved when the budget is unlimited (see Section 5). With a finite budget, we can adjust the value of $C_0$. Recall that $C_0$ represents our tolerance to the sample noise variance of the point estimates and that $h(k)=m$ only if $\widehat{\text{var}}(\mathcal{Y}_m(\cdot))\leq C_0$ at some points. Smaller $C_0$ are more conservative and require a larger number of replications to drive $\widehat{\text{var}}(\mathcal{Y}_m(\cdot))$ down below $C_0$. As a result, the search approaches the objective level slowly. Therefore, to speed up the increase of $h(k)$ to $m$, we can enlarge $C_0$ to ensure that at some design points, $\widehat{\text{var}}(\mathcal{Y}_m(\cdot))\leq C_0$ before the budget runs out. More formally, we can update $C_0$ in the $k$th iteration, denoted by $C_0^k$, as follows:
$$C_0^k=\max\{C_0^{k-1}, \ \ \frac{\widehat{\epsilon}(\widehat{x}_k)N_k(\widehat{x}_k)}{N_k(\widehat{x}_k)+ \frac{A}{|D_k|+\frac{A}{B_k}}} \},  $$
where $\widehat{x}_k$ is the current best, $N_k(\widehat{x}_k) $ is the current number of replications at $\widehat{x}_k$, and $\widehat{\epsilon}(\widehat{x}_k)$ is the noise estimate of $\mathcal{Y}_m(\widehat{x}_k)$. Suppose the number of iterations spent in the next few iterations have the same magnitude as $B_k$. The quantity $\frac{A}{B_k}$ represents an estimate of total possible remaining iterations, and thus the number of new design points selected in the future. It follows that $|D_k|+\frac{A}{B_k}$ is an estimate of the total number of design points at the end of the search. We further assume the remaining replications are evenly distributed and then at the end, $\widehat{x}_k$ can receive $\frac{A}{|D_k|+\frac{A}{B_k}}$ more replications. Therefore,
$\frac{\widehat{\epsilon}(\widehat{x}_k)N_k(\widehat{x}_k)}{N_k(\widehat{x}_k)+ \frac{A}{|D_k|+\frac{A}{B_k}}}$ can be treated as the sample noise variance of the $\alpha_m$-quantile estimator at $x_0^*$. By updating $C_0$ in this way, we aim to make the noise of the $\alpha_m$-quantile estimator at $\widehat{x}_k$ smaller than $C_0$ in the end. As a result, the design set for $\alpha_m$ is non-empty and we reach the objective level. A simpler effort-based rule can be used to approach the $\alpha_m$-th level as well. In this rule, users can define the maximum number of replications that can be spent at the quantiles lower than $\alpha_m$. The algorithm will be forced to go to the $\alpha_m$th level after this budget is exhausted.

\section{eTSSO-QML Convergence Analysis}

This section demonstrates the consistency of the eTSSO-QML algorithm. We first introduce our main assumptions.
\begin{assumption}
	
	(i) There exist $M$ and $f^*$, such that $|v_{\alpha_m}(L(x))|<M$ and $f_x(v_{\alpha_m}(L(x)))>f^*>0$ for all $x\in \mathcal{X}$, where $f_x$ is the probability density function for $L(x)$.
	
	(ii) The variance parameter $\sigma^2$ and the sensitivity parameter $\theta$ of the GP model are bounded away from zero. 
	
	(iii) The model input $R_\epsilon$ is known.
\end{assumption}
Similar assumptions are used in the analysis of the EGO algorithm \citep{jones1998efficient}. Assumption 5.1(i) bounds the optimization objective function and the noise variance of the point estimates. The uniform bound on $f_x(v_{\alpha_m}(L(x)))$ helps ensure the consistency of eTSSO-QML. Assumption 5.1(ii) ensures that the GP model is efficient and reasonable to use. Otherwise, with zero $\sigma^2$, the uncertainties at unsampled inputs would become zero meaning all unobserved points are actually known. With zero $\theta$, the correlation between responses at any two inputs would become zero. In other words, point estimates at the design points would not help predict the responses at unobserved points. Hence, any spatial metamodel would be ineffective under either condition. Assumption 5.1(iii) is used in the convergence proof of the original eTSSO algorithm. This convergence, however, can be affected by the quality of the $R_\epsilon$ estimator. To the best of our knowledge, when $R_\epsilon$ is estimated, convergence has only been studied empirically \citep{kleijnen2012expected,Pedrielli2018etsso} (in the deterministic setting where $R_\epsilon=0$, the convergence has been studied theoretically by \cite{bull2011convergence}). The empirical convergence with estimated $R_\epsilon$ can be partially seen from the numerical tests in Section 6.

The convergence proof for the eTSSO-QML algorithm consists of three parts. First, in Lemma 5.2, we prove that as the iteration increases, the adopted co-kriging model tends to a single-level model for the objective level $v_{\alpha_m}(L(x))$. This is intuitive since each selected design point will be allocated an infinite number of replications as the iteration number increases to infinity. In this case, the model at the objective level is accurate enough so that we may optimize it without leveraging the lower levels.
\begin{lem}
	Under Assumptions 4.1 and Assumption 5.1, there exists $K$ such that for iterations $k>K$, eTSSO-QML reaches a single-level stochastic GP model for $v_{\alpha_m}(L(x))$, i.e., $\pi_k=\{\alpha_m \}$. 
\end{lem}
Second, in Lemma 5.3, we prove that the design points selected by the algorithm are dense in the design space. Generally, convergence proofs for global optimization algorithms require dense design points \citep{torn1989global}.
\begin{lem}
	Under Assumption 4.1, the sequence of design points $D_k$ selected by eTSSO-QML is dense in the design space as $k\rightarrow \infty$. 
\end{lem}
Finally, in Theorem 5.4, we prove the convergence of the overall eTSSO-QML algorithm when the replications at each design point tend to infinity. Recall that the algorithm reports $\widehat{x}_k=\arg \min_{x\in D_k} \mathcal{Y}_m (x) $ as the optimal solution and $\widehat{\mathcal{Y}}_k=\mathcal{Y}_m (\widehat{x}_k) $ as the optimal value found by iteration $k$. Under Assumption 4.1, the number of replications at every selected design point increases to infinity uniformly. Hence, we obtain the following convergence result for eTSSO-QML. 

\begin{thm}
	Under Assumption 4.1, the optimal value found by eTSSO-QML converges to the true global optimum: $\widehat{\mathcal{Y}}_k \rightarrow v_{\alpha_m}(L(x^*))$ w.p.1 as $k\rightarrow \infty$, where $x^*=\arg \min_{x\in \mathcal{X}} v_{\alpha_m}(L(x))$ is the true optimal solution.
\end{thm}

The proof of Theorem 5.4 first considers known hyperparameters $\theta$ and $\sigma^2$. Then we extend this argument to the case where the hyperparameters are estimated. A similar result is proved by \cite{bull2011convergence} for the deterministic setting when the hyperparameters are bounded.

\section{Numerical Experiments}
In this section, we run a few numerical experiments to test the performance of the eTSSO-QML algorithm. In Section 6.1, two simple one-dimensional tests are first presented to illustrate the evolution of the algorithm as designed and in section 6.2, more complicated examples are carried out to further compare eTSSO-QML algorithm with eTSSO-Q algorithm, which directly optimizes the objective quantile function.

\subsection{Two Simple One-Dimensional Illustrating Examples}
As described in Section 4, eTSSO-QML is designed to first search some lower quantile with a small budget to identify the promising regions. After which the algorithm expends more budget into these promising regionsto find the optimum for the higher objective quantile. This section provides two simple examples to examine this design behavior. For a simple and clear illustration of the algorithm, we only consider two quantile levels: $\alpha_1=0.6$, $\alpha_2=0.95$.  

\subsubsection{Experiment 1}
The design space for this example is $\mathcal{X}=[0,1]$. At each $x\in\mathcal{X}$, the loss $L(x)$ is assumed to be normal distributed with mean $m(x)$ and variance $v(x)$:
\[m(x)=5(0.2(x-0.02)+1)\cos(13(x-0.02)) , \ \ v(x)=5x. \]

\begin{figure}[h]	
	\centering
	\begin{minipage}[l]{0.5\textwidth}
		\label{E1}
		\centering
		\includegraphics[width=8cm]{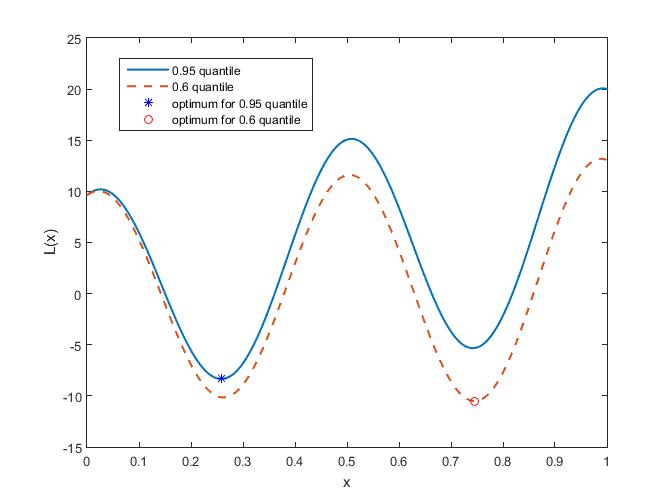}	
		\caption{True quantile functions for Experiment 1}	
	\end{minipage}%
	\begin{minipage}[l]{0.5\textwidth}
		\label{E2}
		\centering
		\includegraphics[width=8cm]{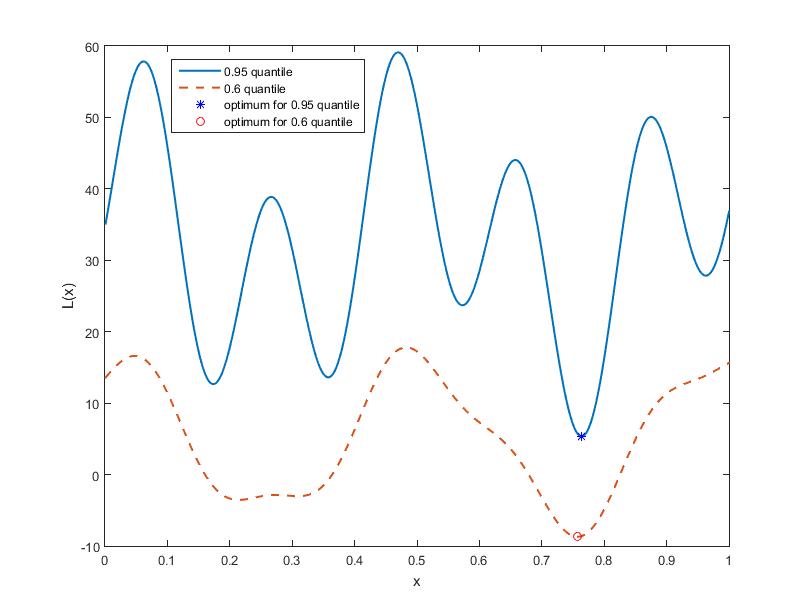}	
		\caption{True quantile functions for Experiment 2}	
	\end{minipage}
\end{figure}

Figure 2 shows the true 0.6 and 0.95 quantile functions as well as their optimums. For this problem, our initial design consists of 6 points selected by Latin Hypercube design, with $r_0=50$, $T=1000$. To show if the algorithm evolves as expected, we provide details about the selected designs and replications assigned to them iteratively (in Table 3). 
\begin{table}[h] \scriptsize
	\centering
	\caption{Design points selected and the number of replications assigned to them in each iteration}
	\begin{tabular}{c|c|c|c|c|c|c|c|c|c|c|c|c}
		\hline
		\multirow{2}{*}{Iteration} & \multicolumn{12}{c}{Design points selected} \\
		\cline{2-13}
		 & 0.085 & 0.2008 &0.3923 &0.5924&0.7057&0.9689&0.721&0.259&0.737&0.748&0.76&0.264\\
		 \hline
		Initial   & 50 & 50 &50 &50&50&50&0&0&0&0&0&0\\
		1 & 0&0&0&0&0&0&50&0&0&0&0&0\\
		2 & 0&0&0&0&0&0&0&50&0&0&0&0\\
		3 & 0&0&0&0&1&0&14&17&70&0&0&0\\
		4 & 1&1&0&0&4&0&7&5&11&70&0&0\\
		5 & 1&1&0&0&5&0&6&5&10&12&64&0\\
		5 & 0&5&1&1&18&1&24&26&30&37&49&110\\
		\hline
	\end{tabular}
\end{table}

In the first 2 iterations, we search the lower quantile and then quickly concentrate more comprehensively into the promising regions (0.2, 0.3) \& (0.7, 0.8) to search for the optimum for the 0.95 quantile. The final optimum found is 0.259 (the true optimum is 0.258). With this example, we can see that the promising regions are correctly identified by first searching the lower level and when we shift the search to the higher quantile, we correctly focus on these promising regions to find the optimum. Another observation is that at the points sampled in the non-promising regions, such as 0.085, 0.3923, 0.5924, the algorithm almost assigns only $r_0$ replications to them. Here we see that the algorithm quickly identifies and eliminates the non-promising regions with a lower quantile model. 

An alternative to this multi-level metamodel search is to directly search the 0.95 quantile with a single model (eTSSO-Q). However, it may be difficult to determine the non-promising regions quickly with this approach as the point estimates of the 0.95 quantile with a similar number of initial runs can be noisy and thus the metamodel built can be misleading. To further investigate this, we conduct another experiment in section 6.1.2.

\subsubsection{Experiment 2}
To illustrate the benefit we can get from eTSSO-QML, we conduct a numerical experiment to compare it with the eTSSO-Q based on a single-level quantile model for the target $\alpha_m$-quantile level. Without too much modification (set $m=h(k)=1$), eTSSO-QML can be easily adapted to eTSSO-Q. 

The mean for $L(x)$ used here is similar with Experiment 1 with variance:
\[ v(x)=10(2+\sin(10\pi x-0.5)) .\]
The true quantile functions are shown in Figure 3. In this experiment, we set $r_0=20$, $T_0=1000$. To mitigate the stochastic nature of the problem, all experiments are conducted with 100 macro-replications. To compare the two algorithms, we further define the true selection as: $|x^*-x_0|<0.035, $ where $x^*$ is the optimum found and $x_0$ represents its true value (0.765). The experiment results for 100 macro-replication are summarized in Table 4.
\begin{table}[h] \small
	\centering
	\caption{Comparison of the two methods}
	\begin{tabular}{lcc}
		\hline
		& eTSSO-QML & eTSSO-Q\\
		Frequency of true selection & 91 & 70 \\
		Average prediction error after initial design & 7.373 (0.6 quantile) & 17.98 (0.95 quantile)\\
		\hline
	\end{tabular}
\end{table}
The results show that eTSSO-QML is much better than eTSSO-Q in terms of true selection. Table 4 further provides the prediction error of the metamodel (measured at another 1000 unsampled points) used by the two algorithms with the initial design. Note that in eTSSO-QML, this initial metamodel is for the 0.6 quantile function while in eTSSO-Q, the metamodel is for the 0.95 quantile function as it directly optimizes this objective level. Seen from the results, with the initial budget, fitting a metamodel for the higher quantile is more inaccurate. In other words, in the beginning, eTSSO-QML utilizes a more accurate surface (0.6 quantile) compared with eTSSO-Q, which uses a very inaccurate 0.95 surface. This inaccurate surface can mislead the search and waste some of the budget on unpromising regions. To show this, we here provide the design inputs that have been selected by the two approaches in one macro-replication run (see Table 5). In this run, in the first four iterations, eTSSO-QML searches 0.6-quantile level ($h(k)=1$) and then goes to 0.95-quantile level with a focus on the promising regions around 0.75. eTSSO-Q, however, seems to still focus on providing a space-filling design (likely due to the poor model estimation throughout) and does not end up near the optimum.

\begin{table}[h] \small
	\centering
	\caption{Comparison of the design inputs selected by the two methods in one run}
	\begin{tabular}{lccccccccc}
		\hline
		eTSSO-QML &0.674 & 0.295& 0.111&0.895&0.987&0.709&0.731&0.732&0.74\\
		eTSSO-QM &0.7&0.5&0.3&0.166&0.108&0.873&0.151\\		
		\hline
	\end{tabular}
\end{table}

With the two simple illustrating examples in this section, we observe that the proposed eTSSO-QML first searches on a lower and more accurate quantile function with a limited budget, and then goes up to the objective level focusing on the promising regions identified. The eTSSO-Q that directly searches the objective quantile function, in contrast, can face a very inaccurate response surface, especially in the beginning, which can then mislead the searching process, resulting in a much more inefficient usage of the budget. 

\subsection{Numerical Tests}

This section employs several more complicated test functions to compare eTSSO-QML with eTSSO-Q. Specifically, here we test if the two algorithms can converge to the global optimums of the test function and how fast they converge. The test functions $F_1$ to $F_3$ used are (in the $d$-dimensional input space):
\[\text{Ackley:}\ \ F_1(x)=-20\exp(-0.2\sqrt{\frac{1}{d}\sum_{i=1}^{d}x_i^2})- \exp(\frac{1}{d}\sum_{i=1}^{d}\cos(2\pi x_i) )+20+\exp(1). \]
\[\text{Rastrigin:} \ \ F_2(x)=10d+\sum_{i=1}^{d}[x_i^2-10\cos(2\pi x_i) ] . \]
\[
\begin{split}
\text{Levy :}\ \ &F_3(x)=\sin^2(\pi \omega_1)+\sum_{i=1}^{d-1}(\omega_i-1)^2[1+10\sin^2(\pi \omega_i+1) ]+ (\omega_d-1)^2[1+\sin^2(2\pi \omega_d)],  \\ &\text{where,} \ \ \omega_i=1+\frac{x_i-1}{4}, \forall i=1,...,d.\end{split}\]
These three functions are all commonly-used test functions for optimization problems. $F_1$ and $F_2$ both have a large number of local optima. $F_3$ is badly-scaled as well as multimodal. Based on these functions, we construct the loss functions $L_1$ to $L_3$ as follows:
\[L_i(x)=F_i+ \text{Lognormal }(0, (1.6+0.01\sum_{i=1}^{d}(x_i-1)^2 )^2 ), \ \ i=1,2, \]
\[L_3(x)=F_3+ \text{Lognormal }(0, (1.6+0.01\sum_{i=1}^{d}x_i^2 )^2 ) . \]
Different with the experiments in section 6.1 that use the normal noise, in this section, we consider log-normal noises to construct the loss functions, since the log-normal distribution is heavier-tailed and thus, increasing the difficulty in estimating the high quantiles. The selected log-normal noises ensure that the global minimizer of the loss function gradually shifts from $[0,...,0]$ to $[1,...,1]$ as the quantile levels increases for $L_1$ and $L_2$. For $L_3$, the minimizer shifts from $[1,...,1]$ to $[0,...,0]$.

The objective level of quantile considered in these examples are 0.99 and starting from 0.6 quantile function in the multi-level algorithm. We select 0.65, 0.7, 0.75, 0.8, 0.85, 0.9, 0.95 as the inter-levels. For each test loss function, both the multi and single level algorithms are run 40 times in the 5-dimensional input region $[-10,10]^5$. Similar to the work of the SKO \citep{huang2006global}, we document the true quantile values at the objective level corresponding to the current found best input $\widehat{x}_k$ in each iteration $k$. The averaged results over these 40 runs are presented in Figure 4. 
\begin{figure}
	\label{average performance}
	\centering
	\begin{minipage}[c]{0.5\textwidth}
		\centering
		\includegraphics[width=8cm]{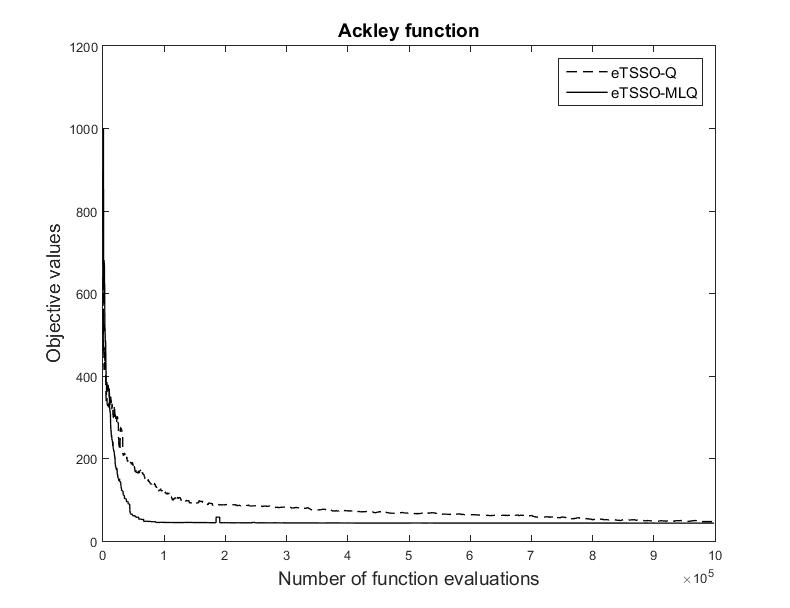}
	\end{minipage}%
	\begin{minipage}[c]{0.5\textwidth}
		\centering
		\includegraphics[width=8cm]{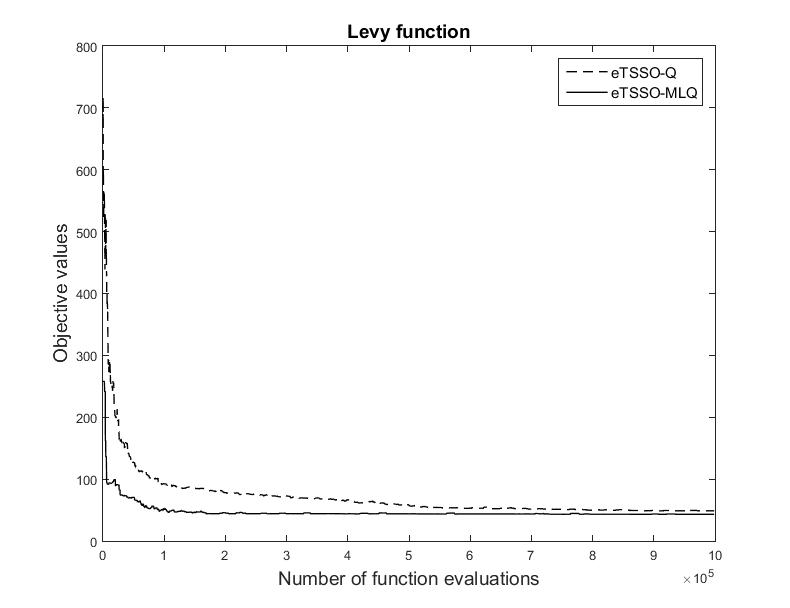}
	\end{minipage}
	\begin{minipage}[c]{0.5\textwidth}
		\centering
		\includegraphics[width=8cm]{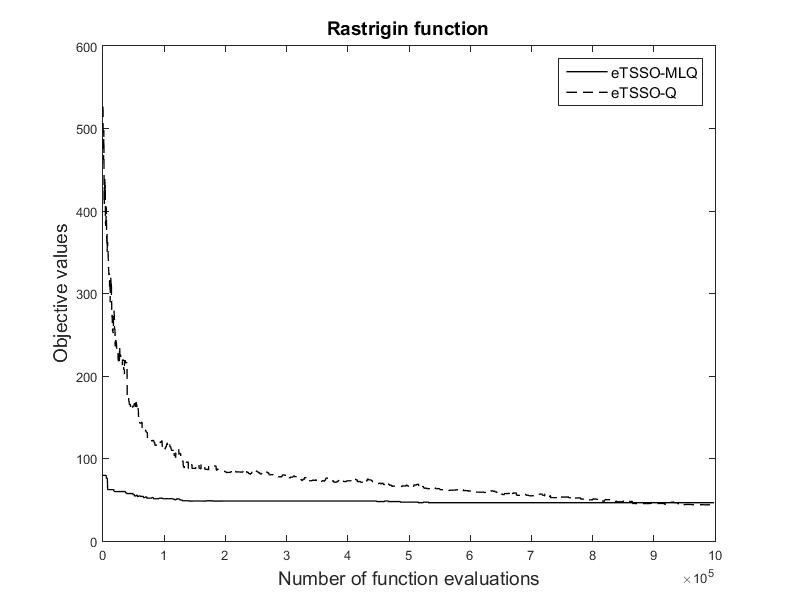}
	\end{minipage}
	\caption{Average performance of the two algorithms}
\end{figure}
We observe from this figure that the two algorithms can converge to the optimal and eTSSO-QML converges faster. As a result, eTSSO-QML often finds the optimal solution with less budget, which is very attractive for expensive simulations.

To compare empirical convergence more clearly, we evaluate the efficiency of the two algorithms through the number of function evaluations cost to find an optimal within a certain relative distance with the true global optimal. Following \citet{barton1984minimization} and \citet{huang2006global}, define $G_k$ in iteration $k$ as follows:
$$G_k=\frac{v_{\alpha_m}(L(x_0))- v_{\alpha_m}(L(\widehat{x}_k)) }{v_{\alpha_m}(L(x_0))- v_{\alpha_m}(L(x^*)) }  ,$$
where $x_0$ is the initial design input of the algorithm and $x^*$ is the true global optimal solution to the objective level quantile function. In this sense, $G_k$ represents the reduction of the gap between the starting value and the current found best over the gap between the starting value and the true global optimum. A larger value of $G_k$ indicates that $\widehat{x}_k$ is closer to $x^*$. As the two algorithms employ multiple initial design inputs and thus multiple starting values of $v_{\alpha_m}(L(x))$, we choose the minimum of these values as $v_{\alpha_m}(L(x_0))$. Similar to the previous works, we use $S_{0.99}$, which is the number of function evaluations until $G_k\geq0.99$ to evaluate the algorithm efficiency. Table 6 lists the percentage of runs reaching $G_{10^6}\geq0.99$. For those runs reaching $G_{10^6}\geq0.99$, we further provide the average of $S_{0.99}$. 
\begin{table}[h] \small
	\centering
	\caption{Percentage of runs reaching $G_{10^6}\geq0.99$ and average of $S_{0.99}$}
	\begin{tabular}{ccc}
		\hline
		Loss function & eTSSO-QML & eTSSO-Q \\
		\hline
		\multirow{2}{*}{$L_1$} & $100\%$ &$ 85\%$ \\
		& 246823.7 &743256.9\\
		\hline
		\multirow{2}{*}{$L_2$} & $100\%$ &$ 95\%$ \\
		& 137913.6 &713404.4\\ 
		\hline
		\multirow{2}{*}{$L_3$} & $97.5\%$ &$ 75\%$ \\
		& 93634.73 &511465.3\\
		\hline
	\end{tabular}
\end{table}
From Table 6, we observe that the eTSSO-QML is more likely and takes less number of evaluations to reach $G\geq 0.99$ and thus it is more efficient than eTSSO-Q. In these tests, we found that eTSSO-QML spends about $10^5$ (1/10 of the total budget) searching the lower quantile functions. With the help of the more accurate and informative lower quantile functions, we can quickly narrow down the searching area and possibly converge faster to the global optimal. The observations are similar with those from the simple examples we obtain from Section 6.1. Moreover, with these examples in Section 6.2, we observe that eTSSO-QML converges with a finite budget, which shows the its convergence empirically in addition to the asymptotic convergence results in Section 5.

\section{Conclusion}
In this paper, we propose eTSSO-QML, a multi-level metamodel based algorithm, to optimize the quantile functions of loss distributions. This algorithm first optimizes lower and informative quantile functions instead of the objective level directly. Compared with higher quantile functions, the lower ones are typically more accurate to estimate, and thus can be easier to optimize with a limited budget. By optimizing the lower quantiles first, we can quickly narrow down the search area to promising regions. As the algorithm proceeds, the quantile level being optimized increases to the objective level and the search process focuses on the promising regions identified by optimizing the previous levels. To achieve this, we first generalize the stochastic co-kriging model to build the metamodel for a series of quantile functions and propose a PMLE approach to prevent the crossing. In the optimization algorithm, we borrow the two-stage framework from the eTSSO algorithm which balances between selecting design inputs and allocating computing budget to them. After integrating the generalized co-kriging metamodel into the algorithm, we always optimize the current highest level of the quantile functions, which increases as the algorithm proceeds and eventually increases to the objective level. Through our numerical tests, we see that the proposed algorithm finds the optimum faster than directly optimizing the objective level and improves the algorithm efficiency.

Optimizing mean functions has been widely studied in the simulation optimization literature. This work demonstrates a possible extension of these optimization algorithms to quantile functions by incorporating quantile estimation techniques. These extensions, like eTSSO-QML, may inherit the advantages and some nice properties from existing algorithms. Furthermore, with the `multi-level' idea, eTSSO-QML can hopefully act as an alternative approach to quantile optimization problems, especially those involving tail quantiles concerning large losses as in finance.

Within our current framework, there are several directions that are worth further investigation. First, as mentioned before, we can consider some more sophisticated approaches, such as the lognormal kriging model \citep{dowd1982lognormal}, to build the metamodel without crossing. Second, other types of selection and allocation rules (apart from the EI criterion and the OCBA technique) can be explored as well. In addition, the essential idea of this work is to leverage some easy-to-get and accurate information when doing optimization. We believe that this idea can be used for optimizing other expensive functions with limited budget, like the Conditional Value-at-Risk and more general families of risk measures. The evaluations of these risk measures at lower risk levels can hopefully provide informative information to help optimize the risk measures at high risk levels, which are typically more expensive, efficiently.

% Acknowledgments here
\section*{Acknowledgements}%
 A preliminary version of this paper was published in the
 Proceedings of the 2018 Winter Simulation Conference. Ng's and Haskell's work was partially supported by Singapore Ministry of Education (MOE) Academic Research Fund (AcRF) Tier 2 grant MOE2015-T2-2-148.
% Leave this (end of acknowledgment)
% References here (outcomment the appropriate case) 
\appendix
\section{Likelihood Function and Estimation of Model Parameters for the Co-Kriging Model}

Given the co-kriging model, the point estimate vector $\mathcal{Y}$ follows a multivariate normal distribution $\mathcal{N}(H\widehat{\beta},R)$. The model parameters can then be straightforwardly estimated by maximizing the loglikelihood function:
\begin{equation*}
\label{likeli}
l(\mathcal{Y},\rho,\theta, \sigma^2)= -\frac{1}{2} \ln(|R|)- \frac{1}{2} {(\mathcal{Y}- H \widehat{\beta})}^TR^{-1}(\mathcal{Y}- H \widehat{\beta}). 
\eqno(A.1)
\end{equation*}  
This approach, however, is to obtain the parameters from models in different levels simultaneously, which involves a multivariate optimization problem. Obviously, this problem becomes more severe as the number of levels increases. To overcome this drawback, we consider more efficient estimation approaches.

As proved by \citet{kennedy2000predicting}, when the observations have no measurement error ($R_\epsilon=0$), the likelihood function of the observation vector can be fully decomposed as follows (to differentiate this case with the stochastic problem (measurement error $R_\epsilon \neq 0$), we use $\mathcal{Z}^T=(\mathcal{Z}_1^T,...,\mathcal{Z}_m ^T)$ to represent the observation vector for the deterministic case where $\mathcal{Z}_1$ is the observation vector for the first level):
$$l(\mathcal{Z},\rho,\theta, \sigma^2)= l_1(\mathcal{Z}_1,\theta_1,\sigma_1^2)+\sum_{j=2}^{m}l_m(\mathcal{Z}_j-\rho_{j-1}\mathcal{Z}_{j-1},\theta_j, \sigma_j^2)  , $$
where $\mathcal{Z}_{j-1}$ consists of the observations of the inputs in $D$ at level ${j-1}$. The vector $\mathcal{Z}_j-\rho_{j-1}\mathcal{Z}_{j-1}$ can be shown to follow a multi-normal distribution $\mathcal{N}(F_j\widehat{\beta}_j,R_j)$, where $F_j={\bf f}_l(D)$, $R_j=A_j$, $\widehat{\beta}_j={ (F_j^TR_j^{-1}F_j )}^{-1}F_j^TR_j^{-1}( \mathcal{Z}_j-\rho_{j-1}\mathcal{Z}_{j-1})$. The function $l_m$ is the loglikelihood for $\mathcal{N}(F_j\widehat{\beta}_j,R_j)$. This decomposition makes optimizing a large scale function $l$ equivalent to optimizing a series of sub problems ($l_1,...,l_m$) with fewer parameters in each, and thus greatly reduces the complexity.

This decomposition, however, is not so straightforward to generalize to the stochastic case. When the noise variance of the observations $R_\epsilon$ is small, which can be accomplished by increasing simulation replications, the decomposition can serve as an approximation of ($A.1$) by ignoring the higher-order terms of $R_\epsilon$ (see the proof in Appendix B). Although this decomposition approach is not as accurate as the standard approach, which optimizes ($A.1$) directly, it greatly reduces the complexity. In practice, optimizing ($A.1$) directly with all model parameters is more difficult and likely to be trapped in sub-optimal regions. A more practical way is to first use the decomposition approach and then treat the optimums found as starting points to apply the standard approach \citep{fricker2013multivariate}.

\section{Proof of Decomposition of loglikelihood}
Here we prove the approximation of the loglikelihood in a stochastic co-kriging model. For simplicity, we only consider two levels and assume $D_1=D_2=D$, $f_1(x)=f_2(x)=0$, which can be easily generalized to more complicated multi-level cases. 
In this simple case, 
$$R=R_z+R_\epsilon=\left(\begin{array}{cc} 
A_1(D) & \rho_1A_1(D) \\ \rho_1A_1(D)^T  & \rho_1^2A_1(D)+A_2(D)
\end{array} \right) + \left(\begin{array}{cc} 
N_1 & N_2 \\ N_2  & N_3
\end{array} \right) ,$$
where $N_1$, $N_3$ represent the noise variance matrix for $\mathcal{Y}_1$ and $\mathcal{Y}_2$, respectively, $N_2$ represents the noise covariance matrix for $\mathcal{Y}_1$ and $\mathcal{Y}_2$.

It can be computed:
\[
\begin{split}
l(\mathcal{Y}_1,\mathcal{Y}_2,\rho_1,\theta,\sigma^2)&=-\frac{1}{2}\ln(|R|)-\frac{1}{2}\mathcal{Y}^TR^{-1}\mathcal{Y}\\
&= -\frac{1}{2}\ln(|A_1+N_1|)-\frac{1}{2}\mathcal{Y}_1^T{(A_1+N_1)}^{-1}\mathcal{Y}_1 -\frac{1}{2}\ln(|F|)\\ &-\frac{1}{2}(\mathcal{Y}_2-(\rho_1A_1^T+N_2)(A_1+N_1)^{-1}\mathcal{Y}_1)^T{F}^{-1}(\mathcal{Y}_2-(\rho_1A_1^T+N_2)(A_1+N_1)^{-1}\mathcal{Y}_1),
\end{split}
\]
where $F:=\rho_1^2A_1+A_2+N_3-{(\rho_1A_1^T+N_2)}{(A_1+N_1)}^{-1}{(\rho_1A_1+N_2)}$.

Suppose that the number of replications at all design points has order $\mathcal{O}(n)$, we find that (see proof of Equation (191), Page 21 from \cite{petersen2012matrix})
$$(A_1+N_1)^{-1} \approx A_1^{-1}- A_1^{-1}N_1A_1^{-1}.$$
It follows that,
$$F\approx\tilde{F}:= A_2+N_3+\rho_1^2N_1-2\rho_1N_2 , $$
$$\mathcal{Y}_2-(\rho_1A_1+N_2)(A_1+N_1)^{-1}\mathcal{Y}_1\approx \mathcal{Y}_2-\rho_1\mathcal{Y}_1.  $$
Therefore, suppose the number of replications at design points has order $\mathcal{O}(n)$, we have
$$l(\mathcal{Y}_1,\mathcal{Y}_2,\rho_1,\theta,\sigma^2)=l_1(\mathcal{Y}_1,\theta_1,\sigma_1^2)+l_2(\mathcal{Y}_2-\rho_1\mathcal{Y}_1,\theta_2,\sigma_2^2)+\mathcal{O}(1/n) , $$ 
where,
$$l_1(\mathcal{Y}_1,\theta_1,\sigma_1^2)=-\frac{1}{2}\ln(|A_1+N_1|)-\frac{1}{2}\mathcal{Y}_1^T{(A_1+N_1)}^{-1}\mathcal{Y}_1 , $$
$$l_2(\mathcal{Y}_2-\rho_1\mathcal{Y}_1,\theta_2,\sigma_2^2)= -\frac{1}{2}\ln(|\tilde{F}|) -\frac{1}{2}(\mathcal{Y}_2-\rho_1\mathcal{Y}_1)^T{\tilde{F}}^{-1}(\mathcal{Y}_2-\rho_1\mathcal{Y}_1).$$ 
Therefore, when $R_\epsilon$ is small, the decomposition above can serve as an approximation of the likelihood function ($A.1$).

\section{Proof of Theorem 3.2}

Suppose $0<\alpha_1<\alpha_2<1$, we prove the consistency and asymptotic unbiasedness of the proposed sectioning covariance estimator for $\mathcal{Y}_1(x)$ and $\mathcal{Y}_2(x)$ with $n$ simulations at $x$. 

First, we refer to Theorem 2.1 from \citet{lin1980asymptotic} on the asymptotic covariance for $\mathcal{Y}_1(x)$ and $\mathcal{Y}_2(x)$:
$$\lim\limits_{n \rightarrow \infty} n\text{cov}(\mathcal{Y}_1(x),\mathcal{Y}_2(x))=\frac{\alpha_1(1-\alpha_2)}{f(v_{\alpha_1})f(v_{\alpha_2})},$$
where $v_{\alpha_1},v_{\alpha_2}$ are the true quantiles and $f$ is the pdf of the underlying distribution. For simplicity, we define $\gamma:=\frac{\alpha_1(1-\alpha_2)}{f(v_{\alpha_1})f(v_{\alpha_2})}$. 

This proof consists of two parts. In C.1, we prove the consistency and asymptotic unbiasedness of $\widehat{\text{cov}}({\mathcal{Y}_1(x),\mathcal{Y}_2(x)})$. In C.2, we derive its MSE.

\subsection{Consistency and Asymptotic Unbiasedness of $\widehat{\text{cov}}({\mathcal{Y}_1(x),\mathcal{Y}_2(x)})$}
Denote the $n$ simulation results as: ${\bf L}:=\{L(x,\xi_1),...,L(x,\xi_n)\}$. Recall that $n_b$ is the batch size and $n_c$ is the number of results in each batch, and thus $n_b \times n_c=n$. Define ${\bf L}^{(j)}:=\{L(x,\xi_{(j-1)n_c+1}),...,L(x,\xi_{jn_c})\} $ as the $j$th batch of simulation results and $\Psi_i$ as the operator to take the sample $\alpha_i$-quantile: $\Psi_i({\bf L})=L_{\llcorner \alpha_i n \lrcorner}$. For example, $\Psi_i({\bf L}^{(j)})$ represents the sample $\alpha_i$-quantile estimator based on the $j$th batch. According to \citet{bahadur1966note} and \citet{chen2016efficient},
\begin{equation*}
\label{eqn:app1}
\Psi_i({\bf L})=v_{\alpha_i}+\frac{1}{n}\sum_{j=1}^{n}\psi_i(L(x,\xi_j))+R_{i,n}, \ \ \psi_i(x)=\frac{\alpha_i-{\bf 1}_{\{x\leq v_{\alpha_i} \} }}{f(v_{\alpha_i})} , \ \,i=1,2,
\eqno(A.2)
\end{equation*} 
where $R_{i,n}$ is the remainder term with $R_{i,n}=O(n^{-3/4}{(\text{log}\text{log}n)}^{3/4})$.
Denote $R_i^{j}$ as the remainder term in $(A.2)$ for $\Psi_i({\bf L}^{(j)})$, $\varphi_i({\bf L})=\frac{1}{n}\sum_{i=1}^{n}\psi_i(L(x,\xi_i))$, $\varphi_i^j=\varphi_i( {\bf L}^{(j)})=\frac{1}{n_c}\sum_{i=1}^{n_c}\psi_i(L(x,\xi_{(j-1)n_c+i}))$, $\bar{\varphi}_i=\frac{1}{n_b}\sum_{j=1}^{n_b}\varphi_i^{j}$ and $\bar{R}_i=\frac{1}{n_b}\sum_{j=1}^{n_b}R_i^{j}$. 

With these notations, the proposed covariance estimator is:
\[
\begin{split}	
&\widehat{\text{cov}}({\mathcal{Y}_1(x),\mathcal{Y}_2(x)}) = \frac{1}{n_b(n_b-1)}\sum_{j=1}^{n_b}\{(\Psi_1({\bf L}^{(j)})-\Psi_1({\bf L}) )(\Psi_2({\bf L}^{(j)})-\Psi_2({\bf L}) )\} \\ &=\frac{1}{n_b(n_b-1)}\sum_{j=1}^{n_b}\{(\Psi_1({\bf L}^{(j)})-\frac{1}{n_b}\sum_{k=1}^{n_b}\Psi_1({\bf L}^{(k)})+\frac{1}{n_b}\sum_{k=1}^{n_b}\Psi_1({\bf L}^{(k)})-\Psi_1({\bf L}) ) \\ & (\Psi_2({\bf L}^{(j)})-\frac{1}{n_b}\sum_{k=1}^{n_b}\Psi_2({\bf L}^{(k)})+\frac{1}{n_b}\sum_{k=1}^{n_b}\Psi_2({\bf L}^{(k)})-\Psi_2({\bf L}) ) \} \\
&= \sigma_1^2+\sigma_2^2,
\end{split}
\]
where we define
$$ \sigma_1^2:=\frac{1}{n_b(n_b-1)}\sum_{j=1}^{n_b}\{(\Psi_1({\bf L}^{(j)})-\frac{1}{n_b}\sum_{k=1}^{n_b}\Psi_1({\bf L}^{(k)})) (\Psi_2({\bf L}^{(j)})-\frac{1}{n_b}\sum_{k=1}^{n_b}\Psi_2({\bf L}^{(k)})) \}, $$
$$\sigma_2^2 := \frac{1}{(n_b-1)}\sum_{j=1}^{n_b}\{(\Psi_1({\bf L})-\frac{1}{n_b}\sum_{k=1}^{n_b}\Psi_1({\bf L}^{(k)})) (\Psi_2({\bf L})-\frac{1}{n_b}\sum_{k=1}^{n_b}\Psi_2({\bf L}^{(k)}))   \}. $$ 

We next derive the asymptotic properties for $\sigma_1^2$ and $\sigma_2^2$ separately in Sections C.1.1 and C.1.2.
\subsubsection{Asymptotic properties for $\sigma_1^2$.}
Note that
$$\sigma_1^2=\frac{1}{n_b(n_b-1)}\sum_{j=1}^{n_b}\{(\varphi_1^j-\bar{\varphi}_1+R_1^j-\bar{R}_1 )(\varphi_2^j-\bar{\varphi}_2+R_2^j-\bar{R}_2 ) \}. $$
It is easy to obtain
\[
\begin{split}
\mathbb{E}[\frac{n}{n_b(n_b-1)}&\sum_{j=1}^{n_b}\{(\varphi_1^j-\bar{\varphi}_1)(\varphi_2^j-\bar{\varphi}_2)\} ] =\frac{n}{n_b}\text{cov}(\varphi_1^1,\varphi_2^1 ) \\ &= \frac{n}{n_bn_c}\text{cov}(\psi_1(L(x,\xi_1)),\psi_2(L(x,\xi_1)) )= \frac{\alpha_1(1-\alpha_2)}{f(v_{\alpha_1})f(v_{\alpha_2})}.
\end{split} 
\]
Recall $\gamma= \frac{\alpha_1(1-\alpha_2)}{f(v_{\alpha_1})f(v_{\alpha_2})} $. We can show that the value inside the expectation converges to $\gamma$ in probability. Specifically, 
\begin{equation*}
\label{eqn:sigma1}
\frac{n}{n_b(n_b-1)}\sum_{j=1}^{n_b}\{(\varphi_1^j-\bar{\varphi}_1)(\varphi_2^j-\bar{\varphi}_2)\} =\frac{n}{n_b(n_b-1)}(\sum_{j=1}^{n_b}\{\varphi_1^j\varphi_2^j   \}-n_b\bar{\varphi}_1\bar{\varphi}_2).
\eqno(A.3)
\end{equation*}
The first term in ($A.3$) is asymptotically equal to $\frac{1}{n_b}\sum_{j=1}^{n_b}(\sqrt{n_c}\varphi_1^j)(\sqrt{n_c}\varphi_2^j) $, with:
$$\mathbb{E}[\sqrt{n_c}\varphi_1^j\sqrt{n_c}\varphi_2^j ]= \mathbb{E}[\psi_1(L(x,\xi_1))\psi_2(L(x,\xi_1))]=\gamma .$$
Therefore, for all $\epsilon>0$, 
\[
\begin{split}
\mathbb{P}(|\frac{1}{n_b}\sum_{j=1}^{n_b}(\sqrt{n_c}\varphi_1^j)(\sqrt{n_c}\varphi_2^j)-\gamma|>\epsilon ) &\leq \frac{1}{n_b}\frac{1}{\epsilon^2} \text{var}[\sqrt{n_c}\varphi_1^j\sqrt{n_c}\varphi_2^j] \\ &\leq \frac{1}{n_b}\frac{1}{\epsilon^2} \mathbb{E}[n_c {\varphi_1^j}^2 n_c {\varphi_2^j}^2 ] \\ &\leq \frac{1}{n_b}\frac{1}{\epsilon^2}\mathbb{E}[n_c^2 {\varphi_1^j}^4 + n_c^2 {\varphi_2^j}^4 ].
\end{split}
\]
According to \citet{chen2016efficient}, $\mathbb{E}[{\varphi_i^j}^4 ]=\mathcal{O}(n_c^{-2})$. The second term in ($A.3$) is asymptotically equivalent to $n_c\bar{\varphi}_1\bar{\varphi}_2$ and
$$n_c\bar{\varphi}_1\bar{\varphi}_2\leq \frac{1}{2}( {(\sqrt{n_c}\bar{\varphi}_1 )}^2+{(\sqrt{n_c}\bar{\varphi}_2 )}^2  ). $$
It is easy to see that,
$$\mathbb{P}[ |\sqrt{n_c}\bar{\varphi}_i| > \epsilon] \leq \frac{n_c \mathbb{E}[\bar{\varphi}_i^2 ]}{\epsilon^2}= \frac{n_c \text{var}[{{\varphi}_i^1} ]}{n_b\epsilon^2}= \frac{ \text{var}[{{\phi}_i} ]}{n_b\epsilon^2}\rightarrow 0, \ \ i=1,2.  $$
It follows that $\frac{n}{n_b(n_b-1)}\sum_{j=1}^{n_b}\{(\varphi_1^j-\bar{\varphi}_1)(\varphi_2^j-\bar{\varphi}_2)\}$ converges to $\gamma$ in probability.

On the other hand, according to the Cauchy-Schwarz inequality, 
\[
\begin{split}
&\mathbb{E}[\frac{n}{n_b(n_b-1)}\sum_{j=1}^{n_b}\{(R_1^j-\bar{R}_1)(R_2^j-\bar{R}_2)  \}  ] \\
&\leq \mathbb{E}[\frac{n}{n_b(n_b-1)} \sqrt{\sum_{j=1}^{n_b}{(R_1^j-\bar{R}_1)}^2}\sqrt{\sum_{j=1}^{n_b}{(R_2^j-\bar{R}_2)}^2}  ]\\
& \leq \sqrt{\frac{n}{n_b(n_b-1)}\mathbb{E}[\sum_{j=1}^{n_b}{(R_1^j-\bar{R}_1)}^2 ]   \frac{n}{n_b(n_b-1)}\mathbb{E}[\sum_{j=1}^{n_b}{(R_2^j-\bar{R}_2)}^2 ] }.
\end{split}
\] 
The second inequality follows the Cauchy-Schwarz inequality applied in probability theory that $|\mathbb{E}[X_1X_2]|^2\leq \mathbb{E}[X_1^2] \mathbb{E}[X_2^2]$, where $X_1$ and $X_2$ are random variables. According to \cite{duttweiler1973mean}, $\frac{n}{n_b(n_b-1)}\mathbb{E}[\sum_{j=1}^{n_b}{(R_1^j-\bar{R}_1)}^2 ]=\mathcal{O}(n_c^{-1/2})$. Therefore, the above expectation converges to zero. Similarly, 
\[ 
\begin{split}
&\frac{n}{n_b(n_b-1)}\sum_{j=1}^{n_b}\{(R_1^j-\bar{R}_1)(R_2^j-\bar{R}_2)  \}  \\
& \leq \sqrt{\frac{n}{n_b(n_b-1)}\sum_{j=1}^{n_b}{(R_1^j-\bar{R}_1)}^2  }\sqrt{\frac{n}{n_b(n_b-1)}\sum_{j=1}^{n_b}{(R_2^j-\bar{R}_2)}^2  }.
\end{split}
\]
According to \cite{chen2016efficient}, $\frac{n}{n_b(n_b-1)}\sum_{j=1}^{n_b}{(R_1^j-\bar{R}_1)}^2 =\mathcal{O}(n_c^{-1/2}(\log \log n_c)^{3/2} )$. It follows that the above quantity converges to zero. 

By using the Cauchy-Schwarz to the cross product terms, $\sigma_1^2$ is shown to be asymptotically unbiased and converges to $\gamma$. 

\subsubsection{Asymptotic properties for $\sigma_2^2$}We next prove the property for $\sigma_2^2$.

Through simple computation,
\[ 
\begin{split}
\sigma_2^2=\frac{1}{n_b^2(n_b-1)}\sum_{j=1}^{n_b}\{R_1^j-R_{1,n} \}\sum_{j=1}^{n_b}\{R_2^j-R_{2,n} \}.
\end{split}
\]
According to \citet{chen2016efficient}, $\frac{1}{n_b^2(n_b-1)}\{\sum_{j=1}^{n_b}(R_1^j-R_{1,n}) \}^2 \rightarrow 0$, $\mathbb{E}[\frac{1}{n_b^2(n_b-1)}\{\sum_{j=1}^{n_b}(R_1^j-R_{1,n}) \}^2]\rightarrow 0$. It can be easily proved that $\sigma_2^2$ converges to 0 in probability and is asymptotically unbiased. 

With Section C.1.1 and C.1.2, following the asymptotic properties for $\sigma_1^2$ and $\sigma_2^2$, it is easy to see the convergency and asymptotic unbiasedness of the proposed covariance estimator.

\subsection{MSE of $\widehat{\text{cov}}({\mathcal{Y}_1(x),\mathcal{Y}_2(x)})$}
We next check the MSE of $\widehat{\text{cov}}({\mathcal{Y}_1(x),\mathcal{Y}_2(x)})$. Its bias is easy to see from the proof in Section C.1 and the squared bias has order $o(n^{-2})$. We next only check the variance of $\widehat{\text{cov}}({\mathcal{Y}_1(x),\mathcal{Y}_2(x)})$.

According to the Cauchy-Schwarz inequality,
\[
\begin{split}
&\text{var}(\widehat{\text{cov}}({\mathcal{Y}_1(x),\mathcal{Y}_2(x)})) \\ 
= &\text{var}(\frac{1}{n_b(n_b-1)}\sum_{j=1}^{n_b}\{(\Psi_1({\bf L}^{(j)})-\Psi_1({\bf L}) )(\Psi_2({\bf L}^{(j)})-\Psi_2({\bf L}) )\} )	\\
\leq &\mathbb{E}[\left\{(\frac{1}{n_b(n_b-1)}\sum_{j=1}^{n_b}\{(\Psi_1({\bf L}^{(j)})-\Psi_1({\bf L}) )(\Psi_2({\bf L}^{(j)})-\Psi_2({\bf L}) )\}  \right\} ^2 ]\\
\leq & \mathbb{E}[\{\frac{1}{n_b(n_b-1)}\sum_{j=1}^{n_b}(\Psi_1({\bf L}^{(j)})-\Psi_1({\bf L}))^2 \} \{\frac{1}{n_b(n_b-1)}\sum_{j=1}^{n_b}(\Psi_2({\bf L}^{(j)})-\Psi_2({\bf L}))^2 \}]	\\
\leq & \mathbb{E}[\frac{1}{n_b(n_b-1)}\sum_{j=1}^{n_b}(\Psi_1({\bf L}^{(j)})-\Psi_1({\bf L}))^2]\mathbb{E}[\frac{1}{n_b(n_b-1)}\sum_{j=1}^{n_b}(\Psi_2({\bf L}^{(j)})-\Psi_2({\bf L}))^2]\\
+ & \sqrt{\text{var}[\frac{1}{n_b(n_b-1)}\sum_{j=1}^{n_b}(\Psi_1({\bf L}^{(j)})-\Psi_1({\bf L}))^2  ]  \text{var}[\frac{1}{n_b(n_b-1)}\sum_{j=1}^{n_b}(\Psi_2({\bf L}^{(j)})-\Psi_2({\bf L}))^2  ]        }
\end{split}
\]
The first inequality follows that $\text{var}[X]\leq\mathbb{E}[X^2]$, the second inequality follows Cauchy-Schwarz inequality and the third inequality follows that $\text{cov}[X_1,X_2]=\mathbb{E}[X_1X_2]-\mathbb{E}[X_1]\mathbb{E}[X_2]\leq \sqrt{\text{var}[X_1]\text{var}[X_2]}$, where $X,X_1,X_2$ are random variables.
From Sections C.1.1 and C.1.2, $\mathbb{E}[\frac{1}{n_b(n_b-1)}\sum_{j=1}^{n_b}(\Psi_i({\bf L}^{(j)})-\Psi_i({\bf L}))^2  ]=o(n^{-1}) $. According to \cite{chen2016efficient}, $\text{var}[\frac{1}{n_b(n_b-1)}\sum_{j=1}^{n_b}(\Psi_i({\bf L}^{(j)})-\Psi_i({\bf L}))^2  ] =o(n^{-2})$. It follows that the variance of the proposed estimator has order $o(n^{-2})$.

\section{Proof of Theorem 3.3}
Consider model (6) for the non-negative response function $W$ where $\theta_0$ is the true value of the hyperparameter for this model. Denote $l_n(\theta)$ as the ordinary loglikelihood function and $Q_n(\theta)$ as the penalized likelihood function:
$$Q(\theta)=-l(\theta)+\lambda \kappa(\theta). $$

We first compute the order of the penalty term $\kappa(\theta)$. The main idea is that, as the design points get denser and denser, the difference between any unobserved design point $x$ with its nearest design point becomes smaller and so does the difference between the predictive value at $x$ and the positive observation at its nearest design point. In this case, the predictive value at $x$ becomes more likely to be positive (as it becomes more and more close to a positive value). For simplicity, let the design space be one-dimensional. Nonetheless, this proof can be easily generalized to multi-dimensional case. 

For each unobserved $x\in \mathcal{X}$, let $x_0\in \{D\}$ be the nearest design point to $x$ (if there are more than one nearest point, pick any one):
$$x_0:=\arg\min_{x'\in D} |x'-x|.  $$
Further denote $h_0$ as the maximum of the distance between any unexplored input with its nearest design input:
$$h_0:=\sup_{x\in \mathcal{X}\setminus D} \ \  \inf_{x^0 \in  D }|x-x^0  |. $$
Within a fixed domain, as the design points become dense, $h_0 \rightarrow 0$. In other words, there exists a sequence $c_n$ such that $c_n \rightarrow \infty$ as $n\rightarrow \infty$ and $h_0=\mathcal{O}(c_n^{-1})$.

The predictor for the deterministic GP model considered here has similar form with (2) by setting noise variance matrix $R_\epsilon=0$ and the number of levels as 1. Moreover, the predictive value at $x_0$ is exactly the observation here owing to the interpolation property of the deterministic GP model: 
$$\mathcal{W}(x_0)= f(x_0)^T\beta +t(x_0)^T R^{-1}(\mathcal{W}-F\beta),$$
we see that:
$$\widehat{W}(x)=\mathcal{W}(x_0)+(f(x)^T-f(x_0)^T)\beta+(t(x)^T-t(x_0)^T )R^{-1}(\mathcal{W}-F\beta),  $$
where $t(x)$ is the covariance vector between $x$ and the design points. For any entry in $f(x)$ and $t(x)$, considering Taylor expansion, we have:
$$f(x)_1-f(x_0)_1=f'(x_0)|x-x_0|+o(|x-x_0|)=\mathcal{O}(h)=\mathcal{O}(c_n^{-1}), $$
$$t(x)_1-t(x_0)_1=t'(x_0)|x-x_0|+o(|x-x_0|)=\mathcal{O}(h)=\mathcal{O}(c_n^{-1}). $$
Here, we assume that for each $\theta$, the first derivative of $f$ and $t$ is bounded within the design domain, which is valid in most cases. Considering that $f(x)^T\beta+t(x)^TR^{-1}(\mathcal{Z}-F\beta)$ is of order $\mathcal{O}(1)$, we see that $\widehat{W}(x)-\mathcal{W}(x_0)=\mathcal{O}(c_n^{-1})$. In this case, the penalty term $\kappa(\theta)=\mathcal{O}(c_n^{-1})$ for every possible $\theta$.

For the ordinary MLE, according to \citet{yi2011penalized}, under the similar regularity conditions, there exists a solution $\widehat{\theta}_n$ to $l(\theta)=0$, which is consistent for $\theta_0$ as $n\rightarrow\infty$. According to \citet{stein2012interpolation} and \citet{li2005analysis}, for this series of $\widehat{\theta}_n$, $||\widehat{\theta}_n-\theta||=\mathcal{O}_p(n^{-1/2}). $

The remaining proof is quite similar to Theorem 1 from \citet{fan2001variable}. We need to show that for all $\epsilon>0$, there exists a large constant $C$, such that:
\begin{equation*}
\label{eqn:Pp}
\mathbb{P}\{\inf\limits_{|| {\bf u} ||=C} Q(\theta_0+n^{-1/2}{\bf u}) > Q(\theta_0)  \} \geq 1-\epsilon.
\eqno(A.4)
\end{equation*}

This shows that there exists a local minimum of $Q$ within the ball $\{\theta_0+n^{-1/2}{\bf u}: || {\bf u} ||\leq C \}$ with probability no less than $1-\epsilon$. It follows that there is a local minimizer of Q satisfying $||\widehat{\theta}_n-\theta_0||=\mathcal{O}_p(n^{-1/2}).  $

Define $D_n({\bf u}):= Q(\theta_0+n^{-1/2}{\bf u})-Q(\theta_0)$, we have,
\[
\begin{split}
D_n({\bf u}&)=-l_n(n^{-1/2}{\bf u})+l(\theta_0)+\lambda(\kappa(n^{-1/2}{\bf u})-\kappa(\theta_0)) \\ 
&= -n^{-1/2}l_n'(\theta_0){\bf u}+\frac{1}{2}{\bf u}^T\mathbb{I}_n(\theta_0){\bf u}n^{-1}(1+\mathcal{O}_p(1))+\lambda\kappa(\sqrt{n}{\bf u})-\lambda \kappa(\theta_0),
\end{split}
\]
where $\mathbb{I}_n(\theta_0)$ is the fisher information matrix. According to \citet{yi2011penalized}, $n^{-1/2}l_n'(\theta_0)=\mathcal{O}_p(1)$, $\mathbb{I}_n(\theta_0)=\mathcal{O}_p(n)$ and thus the second term has order $\mathcal{O}_p(1)$. The third term is positive and the last term has order $\mathcal{O}_p(c_n^{-1})$. By choosing a sufficiently large $C$, the second term dominates the first and the last term. It follows that $D_n({\bf u})>0$ and ($A.4$) holds.

\section{Proof of Lemma 5.2}

At iteration $k$, all selected design points have been allocated at least $r_k>k$ iterations in the proposed algorithm. Consider a design point $x_0\in D_k$. Recall that ($A.2$) states that
$$\mathcal{Y}_m(x_0)=v_{\alpha_m}(L(x_0))+ \frac{1}{N_k(x_0)} \sum_{j=1}^{N_k(x_0) } \psi_m(L(x_0,\xi_j) ) +R_{N_k(x_0)}, \ \ \psi_m(L(x_0,\xi_j) )=\frac{\alpha_m-1_{\{L(x_0,\xi_j)\leq  v_{\alpha_m}(L(x_0))\}}}{f_{x_0}(v_{\alpha_m}(L(x_0)))},$$
where $N_k(x_0)$ is the number of replications assigned to $x_0$ by iteration $k$. It follows that $\mathcal{Y}_m(x_0)$ has variance:

\[\begin{split}
\text{var}(\mathcal{Y}_m(x_0))= \frac{1}{N_k(x_0)}\text{var}(\psi_m(L(x_0,\xi)))+\text{var}(R_{N_k(x_0)})=\frac{\alpha_m(1-\alpha_m)}{N_k(x_0)f^2_{x_0}(v_{\alpha_m}(L(x_0)))}+\text{var}(R_{N_k(x_0)}).
\end{split}\]

According to \cite{duttweiler1973mean}, $\mathbb{E}[R_{N(x_0)}^2]\approx N(x_0)^{-3/2}f_{x_0}^{-2} (v_{\alpha_m}(L(x_0)) ) (2\alpha_m(1-\alpha_m)/\pi)^{1/2}$. Recall that $f_x(v_{\alpha_m}(L(x)))>f^*$ for all $x \in\mathcal{X}$. Under Assumption 4.1, $N(x_0)\geq r_k>k$ for all $x \in D_k$ in iteration $k$. Therefore, for all $x_0\in D_k$,
\begin{equation*}
\label{AppE}
\text{var}(\mathcal{Y}_m(x_0))\leq \frac{\alpha_m(1-\alpha_m)}{N_k(x_0)f^2_x(v_{\alpha_m}(L(x)))}+ \mathbb{E}[R_{N(x_0)}^2] \leq \frac{\alpha_m(1-\alpha_m)}{r_k(f^*)^2}+\frac{(2\alpha_m(1-\alpha_m))^{1/2}}{r_k^{3/2}(f^*)^2\pi^{1/2}  } := P_k.
\eqno(A.5)
\end{equation*}
The first inequality holds since $\text{var}(R_{N_k(x_0)})  \leq \mathbb{E}[R_{N(x_0)}^2]$. We see that $P_k$ does not depend on $x_0$ and furthermore $P_k\rightarrow 0$ as $k\rightarrow \infty$. Note that $P_k$ is an upper bound for $\text{var}(\mathcal{Y}_m(x_0))$ in iteration $k$ for all $x_0\in D_k$. 
It follows that, the noise variance for the $\alpha_m$-quantile estimators at all design points in $D_k$ tends to zero uniformly as $k\rightarrow \infty$. Therefore, there exists a large value $K$ that does not depend on $x$ such that when $k>K$, $\text{var}(\mathcal{Y}_m(x_0))<C_0$ for any design point $x_0\in D_k$ (under Assumption 5.1(iii), we may use the true value of $\text{var}(\mathcal{Y}_m(x_0))$, instead of its estimate (4), to examine the quality of the point estimate). Therefore, in iterations $k>K$ of eTSSO-QML, we only use a single-level model for the objective level. In this case, the model at the objective level is accurate enough such that we may optimize it without leveraging the lower levels.

\section{Proof of Lemma 5.3}

According to Lemma 1, there exists a large number $K$ such that for iterations $k>K$, $\text{var}(\mathcal{Y}_m(x_0))<C_0$ for any design point $x_0\in D_k$ and thus eTSSO-QML adopts a single-level model for the objective level. In this proof, we suppose $k>K$ and omit the subscript $l$ in equations (2) and (7). Denote the EI function in iteration $k$ as ${T}_k(x)$ where:
$$ {T}_k(x)=\widehat{s}_k(x) \phi (\frac{y_k^*-\widehat{\mu}_k(x)}{\widehat{s}_k(x)})+ (y_k^*-\widehat{\mu}_k(x))\Phi(\frac{y_k^*-\widehat{\mu}_k(x)}{\widehat{s}_k(x)}), $$
where $y_k^*$ is the current best objective value. For ease of exposition, we write $\widehat{\mu}_k(x)$ and $\widehat{s}^2_k(x)$ (as there is only one level $\alpha_m$ here, the subscript in $s^2_k(x)$ does not represent the level but the iteration number) as to denote the predictor (2) and predictive variance (7) obtained from the single-level GP model for the objective level.

The proof of this Lemma follows that of Theorem 1 from \cite{locatelli1997bayesian}. It can be divided into three parts. In Section F.1, we find an upper bound for ${T}_k(x)$ at any unobserved point $x\in \mathcal{X} \setminus D_k$. This upper bound depends on the the nearest design point to $x$. Intuitively, if $x$ is very close to a design point $x_0$, the uncertainty at $x$ should be small since its response has a large correlation with $x_0$. As a result, the expected improvement we get from observing the response at $x$ should be small. In Section F.2, we show how to construct a region around any design point where ${T}_k(x)$ is bounded above by a threshold $c$. Finally, in Section F.3, we apply Lemma 1 and Theorem 1 from \cite{locatelli1997bayesian} to prove that the design points are dense.

\subsection{Upper bound for $ {T}_k(x)$}
According to Assumption 5.1(i), the true baseline quantile function is bounded. We may select a large enough value $M$ such that the predictor $\widehat{\mu}_k(x)$ and the true quantile value $v_{\alpha_m}(x)$ are constrained in $(-M,M)$, for all $x\in \mathcal{X}$, for all $k$. 

In iteration $k>K$, consider an unknown point $x$, then we must have:
$$y_k^*-\widehat{\mu}_k(x)< 2M.  $$
Through simple computation of the partial derivatives of $T_k(x)$, we find that $\frac{\partial T_k(x)}{\partial (y_k^*-\widehat{\mu}_k(x))} >0 $ and $\frac{\partial T_k(x)}{\partial \widehat{s}_k(x_0)} >0 $. As $T_k(x)$ is an increasing function of $y_k^*-{\mu}_k(x)$, we have: 
$$T_k(x)\leq \widehat{s}_k(x)\phi(\frac{2M}{\widehat{s}_k(x)})+2M\Phi(\frac{2M}{\widehat{s}_k(x)}) := P_k(x) .$$

Define $x_0:=\arg \max_{x'\in D_k} \text{corr}(x,x')$, which is the design point with the largest correlation with $x$. Denote the covariance matrix in (7) as $R_z=\left[\begin{array}{cc}
R_{11} & R_{12} \\ R_{21}  & R_{22}
\end{array} \right]$, where $R_{11}=\sigma^2$ is the variance of the spatial process at $x_0$, $R_{21}= R_{12}^T$ is the $(|D_k|-1) \times 1$ covariance vector of the spatial process between $x_0$ and the remaining design points $D_k \setminus \{x_0\}$, and $R_{22}$ is the $(|D_k|-1) \times (|D_k|-1)$ covariance matrix for the spatial process at $D_k \setminus \{x_0\}$. Further denote $t(x)=(\sigma^2\text{corr}(x,x_0), t_2^T(x) )^T$, where $t_2(x)$ is the $(|D_k|-1) \times 1$ response covariance vector between $x$ and $D_k \setminus \{x_0\}$. We find that for iteration $k>K$
\[\begin{split}
\widehat{s}^2_k(x) &=\sigma^2- t(x)^TR_z^{-1}t(x)+\zeta(x)^T{(H^TR_z^{-1}H)}^{-1}\zeta(x) \\ &= \sigma^2-\frac{\sigma^4\text{corr}(x,x_0)^2}{R_{11}} - \Gamma^TP\Gamma +\zeta(x)^T{(H^TR_z^{-1}H)}^{-1}\zeta(x)\\
&\leq \sigma^2-\sigma^2 \text{corr}(x,x_0)^2 +\zeta(x)^T{(H^TR_z^{-1}H)}^{-1}\zeta(x),\\
\end{split}
\] 
where $\Gamma=R_{21}R_{11}^{-1}\sigma^2\text{corr}(x,x_0)-t_2(x)$ and $P^{-1}= R_{22}-R_{21}R_{11}^{-1}R_{12}$. The inequality holds since $P^{-1}$ is positive definite (as $P$ is a covariance matrix of the responses at $D_k \setminus \{x_0\}$ given $x_0$, which is symmetric and positive definite). For the last term, we recall that $\zeta(x)=h(x)-t(x)^TR_z^{-1}H$ from (7). For the commonly used constant mean function $h(x)=1$, $t(x)^TR_z^{-1}H$ is the GP prediction at $x$ given observation vector $H$. Following the same procedure as in Appendix D, we have that $h(x_0)-t(x)^TR_z^{-1}H=\mathcal{O}(|x-x_0|)$. As $h(x_0)=h(x)=1$, it follows that $|\zeta(x)|=|h(x)-t(x)^TR_z^{-1}H|=\mathcal{O}(|x-x_0|)$. In this case, we can select a value $M_1$ such that $|\zeta(x)|<M_1|x-x_0|$. Moreover, we can check that ${(H^TR_z^{-1}H)}^{-1}<\sigma^2+C_0$. Therefore $\zeta(x)^T{(H^TR_z^{-1}H)}^{-1}\zeta(x) = {(H^TR_z^{-1}H)}^{-1} |\zeta(x)|^2<M_1^2(\sigma^2+C_0) |x-x_0|^2$. Denote $M_2$ as $M_1^2(\sigma^2+C_0)$, then we have $\zeta(x)^T{(H^TR_z^{-1}H)}^{-1}\zeta(x)< M_2 |x-x_0|^2$. For a general mean function $h$, the proof follows the same reasoning. 

Define $\widehat{s}^2_{k0}(x;x_0):=\sigma^2-\sigma^2  \text{corr}(x,x_0)^2+ M_2 |x-x_0|^2$. We can see that $ \widehat{s}^2_k(x)\leq \widehat{s}^2_{k0}(x;x_0)$. Moreover, as $\text{corr}(x,x_0)$ increases as the distance between $x$ and $x_0$ decreases, we see that $\widehat{s}^2_{k0}(x;x_0)$ increases as the distance between $x$ and $x_0$ increases.

As $P_k(x)$ is an increasing function of $\widehat{s}_k(x)$, we have:
$$T_k(x)\leq P_k(x) \leq \widehat{s}_{k0}(x;x_0)\phi(\frac{2M}{\widehat{s}_{k0}(x;x_0)})+2M\Phi(\frac{2M}{\widehat{s}_{k0}(x;x_0)}) := Q_k(x;x_0).$$ 
Since $\frac{\partial Q_k(x;x_0)}{\partial \widehat{s}_{k0}(x;x_0)}>0$ and the fact that $\widehat{s}^2_{k0}(x;x_0)$ increases as the distance between $x$ and $x_0$ increases, we see that $Q_k(x;x_0)$ increases as the distance between $x$ and $x_0$ increases.

\subsection{Local Region Centered at Design Points with Bounded $ {T}_k(x)$}
Considering any design point $x_0\in D_k$. Based on the bound $Q_k(x;x_0)$, we may construct a region, $R(x_0,c) $, containing $x_0$ defined by:
$$R(x_0,c)=\{x \in \mathcal{X} | Q_k(x;x_0)< c\}.$$ 
From the proof in Appendix F.1, we find that $ Q_k(x;x_0)$ decreases as the distance between $x$ and $x_0$ decreases and that $Q_k(x;x_0)\rightarrow0  $ as $x\rightarrow x_0$. Therefore, for any value of $c$, $R(x_0,c)$ is a region centered at $x_0$, such that $T_k(x)< c$, for all $x\in R(x_0,c)$. 

\subsection{Proof of Density}
To prove the density of the design points, we deploy Lemma 1 and Theorem 1 from \cite{locatelli1997bayesian}. Specifically, we consider the following stopping rule in the algorithm:

-- \emph{Stopping Rule}: The algorithm stops when the maximum of ${T}_k(x)$ is smaller than some pre-defined threshold $c$.

We can see that with this stopping rule, the points in $R(x_0,c) $ will never be selected in iteration $k>K$. From Theorem 1 in \cite{locatelli1997bayesian}, the algorithm will terminate within a finite number of design points for any given $c$. However, we assume infinite budget so that the algorithm does not stop within finitely many iterations. To achieve this, similar to the procedure in \cite{locatelli1997bayesian}, once the algorithm stops, we decrease the value of threshold $c$ to ensure that the maximum of ${T}_k(x)$ is larger that the updated $c$. Thus, the condition of the stopping rule is not met and the algorithm continues. To finish the proof, we directly use the results of Lemma 1 from \cite{locatelli1997bayesian} that the design points will be dense everywhere in $\mathcal{X}$ if the threshold value $c$ keeps decreasing. 

The above result is proved with known parameters $\sigma^2$ and $\theta$. However, when the parameters are estimated, the above proof still goes through under Assumption 5.1(ii). Specifically, the estimated values of these parameters will influence the size of the region $R(x_0,c)$. With bounded values of $\widehat{\sigma}^2$ and $\widehat{\theta}$, the region should be nonempty and the proof will continue to hold.

\section{Proof of Theorem 5.4}

We prove that $\widehat{\mathcal{Y}}_k \rightarrow v_{\alpha_m}(L(x^*))$ w.p.1 as $k\rightarrow \infty$. Equivalently, we prove that $\lim_{n\rightarrow\infty} \mathbb{P}(\cup_{k=n}^\infty \{|\widehat{\mathcal{Y}}_k - v_{\alpha_m}(L(x^*))|>\delta \})=0$, for all $\delta>0$. Recall that $\widehat{\mathcal{Y}}_k=\mathcal{Y}_m(\widehat{x}_k)$, where $\widehat{x}_k=\arg \min_{x\in D_k} \mathcal{Y}_m(x)$ is the observed best point within the design set. Define $x_k^*:=\arg \min_{x\in D_k} v_{\alpha_m}(L(x))$, which is the true best point within the design set. This proof is divided into three parts. In Section G.1, we prove that $\widehat{\mathcal{Y}}_k - v_{\alpha_m}(L(x_k^*))\rightarrow 0$ w.p.1. This is to correctly identify the best points within the design set. In Section G.2, we prove that $v_{\alpha_m}(L(x_k^*)) \rightarrow  v_{\alpha_m}(L(x^*))$, which ensures the true optimum within the design set tends to the true global minimum. In Section G.3, we combine the proofs from G.1 and G.2 to finish the convergence proof.

\subsection{Proof that $\widehat{\mathcal{Y}}_k - v_{\alpha_m}(L(x_k^*))\rightarrow 0$ w.p.1 as $k\rightarrow \infty$}
We prove that $\lim_{n\rightarrow\infty} \mathbb{P}(\cup_{k=n}^\infty \{|\widehat{\mathcal{Y}}_k - v_{\alpha_m}(L(x_k^*))|>\frac{\delta}{2}\})=0$ for all $\delta>0$. To show this, we verify the following sufficient condition $\sum_{k=1}^{\infty}\mathbb{P}[|\widehat{\mathcal{Y}}_k - v_{\alpha_m}(L(x_k^*))>\frac{\delta}{2} | ]<\infty$ (Theorem 7.5, \cite{pishro2016introduction}). For all $\delta>0$,
\begin{equation*}
\label{equationProof3}
\begin{split} 
& \mathbb{P}[|\mathcal{Y}_m(\widehat{x}_k)-  v_{\alpha_m}(L(x_k^*))|>\frac{\delta}{2}]\\  = & \mathbb{P}[|\mathcal{Y}_m(\widehat{x}_k)- v_{\alpha_m}(L(\widehat{x}_k)) + v_{\alpha_m}(L(\widehat{x}_k)) - v_{\alpha_m}(L(x_k^*)) | >\frac{\delta}{2}]\\ < &\mathbb{P}[|\mathcal{Y}_m(\widehat{x}_k)- v_{\alpha_m}(L(\widehat{x}_k))|>\frac{\delta}{4}]+ \mathbb{P}[|v_{\alpha_m}(L(\widehat{x}_k)) - v_{\alpha_m}(L(x_k^*))|>\frac{\delta}{4}]. \ \ \ \ \ (A.6)
\end{split}       
\end{equation*}

We bound the first term in ($A.6$) as follows. Under Assumption 4.1, the accumulated number of replications at each input $x\in D_k$, $N_k(x)>r_k$  and $r_k \rightarrow \infty$ as $k \rightarrow \infty$. According to Bahadur's representation \citep{kiefer1967bahadur}:
$$\mathcal{Y}_m(x)= v_{\alpha_m}(L(x))+\frac{1}{N(x)}\sum_{i=1}^{N(x)}\psi(L(x,\xi_i))+R_{N(x)},$$
where $N(x)$ is the total number of replications at $x$, $\psi(L(x,\xi_i)):=\frac{\alpha_m -\mathbf{1}_{ \{L(x,\xi_i)< v_{\alpha_m}(L(x))\}}}{f_x(v_{\alpha_m}(L(x)))}$ and $R_{N(x)}$ is the remainder term. Therefore, for all $x\in D_k$, we have that 
\[
\begin{split}
\mathbb{P}[|\mathcal{Y}_m(x)- v_{\alpha_m}(L(x))|>\frac{\delta}{4}] &=\mathbb{P}[|\frac{1}{N(x)}\sum_{i=1}^{N(x)}\psi(L(x,\xi_i))+R_{N(x)}|>\frac{\delta}{4}]\\ &<\mathbb{P}[|\frac{1}{N(x)}\sum_{i=1}^{N(x)}\psi(L(x,\xi_i))|>\frac{\delta}{8}]+\mathbb{P}[| R_{N(x)}|>\frac{\delta}{8}] .    
\end{split}
\] 
Through simple computation, we find $\mathbb{E}[ \psi(L(x,\xi_i))]=0$ and $\text{var}[ \psi(L(x,\xi_i))]=\frac{\alpha (1-\alpha)}{f_x^2(v_{\alpha_m}(L(x))) }$, and thus by Chebychev's inequality, we have
$$\mathbb{P}[|\frac{1}{N(x)}\sum_{i=1}^{N(x)}\psi(L_i(x))|>\frac{\delta}{8}]< \frac{64 \alpha (1-\alpha)}{N(x) \delta^2 f_x^2(v_{\alpha_m}(L(x)))} < \frac{64 \alpha (1-\alpha)}{r_k \delta^2 (f^*)^2},$$ (recall $f^*< f_x(v_{\alpha_m}(L(x)))$ for all $x\in \mathcal{X} $ by Assumption 5.1(i)). On the other hand, 
for the remainder $R_{N(x_0)}$, we have
$$\mathbb{P}[| R_{N(x)}|>\frac{\delta}{8}]=\mathbb{P}[| R_{N(x)}|^2>\frac{\delta^2}{64}] \leq \frac{\mathbb{E}[R^2_{N(x)}]}{\delta^2/64}<\frac{64(2\alpha_m(1-\alpha_m))^{1/2}}{\delta^2 r_k^{3/2}(f^*)^2\pi^{1/2}  }.  $$
The first inequality holds by Chebychev's inequality and the second inequality holds by the same reasoning with ($A.5$). It follows that for all $x\in D_k$
$$ \mathbb{P}[|\mathcal{Y}_m(x)- v_{\alpha_m}(L(x))|>\frac{\delta}{4}]<\frac{64 \alpha (1-\alpha)}{r_k \delta^2 (f^*)^2}+\frac{64(2\alpha_m(1-\alpha_m))^{1/2}}{\delta^2 r_k^{3/2}(f^*)^2\pi^{1/2}  } < \frac{64 }{r_k \delta^2 (f^*)^2} \left(\alpha (1-\alpha) +\frac{(2\alpha_m(1-\alpha_m))^{1/2}}{\pi^{1/2}  }  \right).$$ 
Therefore,
\[
\begin{split}
&\mathbb{P}[\max_{x\in D_k}|\mathcal{Y}_m(x)- v_{\alpha_m}(L(x))|>\frac{\delta}{4}] \\ \leq&\sum_{i=1}^{k+|D_0|} \mathbb{P}[|\mathcal{Y}_m(x_i)- v_{\alpha_m}(L(x_i))|>\frac{\delta}{4}] \\ <& \frac{64 (k+|D_0|)}{r_k  \delta^2 (f^*)^2} \left(\alpha (1-\alpha) +\frac{(2\alpha_m(1-\alpha_m))^{1/2}}{\pi^{1/2}  } \right).
\end{split}
\] 
With this inequality, we see that
$$ \mathbb{P}[|\mathcal{Y}_m(\widehat{x}_k)- v_{\alpha_m}(L(\widehat{x}_k))|>\frac{\delta}{4}]\leq\mathbb{P}[\max_{x\in D_k}|\mathcal{Y}_m(x)- v_{\alpha_m}(L(x))|>\frac{\delta}{4}] < \frac{64 (k+|D_0|)}{r_k  \delta^2 (f^*)^2} \left(\alpha (1-\alpha) +\frac{(2\alpha_m(1-\alpha_m))^{1/2}}{\pi^{1/2}  } \right),$$ 
$$ \mathbb{P}[|\mathcal{Y}_m(x_k^*)- v_{\alpha_m}(L(x_k^*))|>\frac{\delta}{4}]\leq\mathbb{P}[\max_{x\in D_k}|\mathcal{Y}_m(x)- v_{\alpha_m}(L(x))|>\frac{\delta}{4}] < \frac{64 (k+|D_0|)}{r_k  \delta^2 (f^*)^2} \left(\alpha (1-\alpha) +\frac{(2\alpha_m(1-\alpha_m))^{1/2}}{\pi^{1/2}  } \right).$$

Now we bound the second term in ($A.6$). Define sets $A_k:=\{|\mathcal{Y}_m(\widehat{x}_k)-v_{\alpha_m}(L(\widehat{x}_k))|\leq\frac{\delta}{9} \}$ and $B_k:=\{|\mathcal{Y}_m(x_k^*)-v_{\alpha_m}(L(x_k^*))|\leq\frac{\delta}{9}\}$ for all $k\geq0$. We note that 
\[
\begin{split}
&\mathbb{P}[|v_{\alpha_m}(L(\widehat{x}_k)) - v_{\alpha_m}(L(x_k^*))|>\frac{\delta}{4}] \\
=&\mathbb{P}[\{|v_{\alpha_m}(L(\widehat{x}_k)) - v_{\alpha_m}(L(x_k^*))|>\frac{\delta}{4}\} \cap \{ A_k\cap B_k\}] +\mathbb{P}[\{|v_{\alpha_m}(L(\widehat{x}_k)) - v_{\alpha_m}(L(x_k^*))|>\frac{\delta}{4}\} \cap \{ A_k\cap B_k\}^\complement].
\end{split}
\] 
We prove that the first term is zero by contradiction. When $v_{\alpha_m}(L(\widehat{x}_k)) - v_{\alpha_m}(L(x_k^*))\geq\frac{\delta}{4}$, as $|\mathcal{Y}_m(\widehat{x}_k)-v_{\alpha_m}(L(\widehat{x}_k))|\leq\frac{\delta}{9}$ (set $A_k$) and $|\mathcal{Y}_m(x_k^*)-v_{\alpha_m}(L(x_k^*))|\leq\frac{\delta}{9}$ (set $B_k$), it must be that $\mathcal{Y}_m(\widehat{x}_k)>\mathcal{Y}_m(x_k^*)$. This inequality contradicts the fact that $\widehat{x}_k$ is the best observed point at iteration $k$, i.e., $\widehat{x}_k=\arg \min_{x\in D_k} \mathcal{Y}_m(x)$. It follows that the first term is 0. For the second term, we see that
\[
\begin{split}
&\mathbb{P}[\{|v_{\alpha_m}(L(\widehat{x}_k)) - v_{\alpha_m}(L(x_k^*))|>\frac{\delta}{4}\} \cap \{ A_k\cap B_k\}^\complement]\\ < &\mathbb{P}[\{ A_k\cap B_k\}^\complement] = 1-\mathbb{P}[ A_k\cap B_k]  <2-\mathbb{P}[ A_k] -\mathbb{P}[ B_k] < \frac{648 (k+|D_0|)}{r_k  \delta^2 (f^*)^2} \left(\alpha (1-\alpha) +\frac{(2\alpha_m(1-\alpha_m))^{1/2}}{\pi^{1/2}  } \right).
\end{split}
\] 
The last inequality follows because $1-\mathbb{P}[A_k]=\mathbb{P}[|\mathcal{Y}_m(\widehat{x}_k)-v_{\alpha_m}(L(\widehat{x}_k))|>\frac{\delta}{9} ]< \frac{324 (k+|D_0|)}{r_k  \delta^2 (f^*)^2} \left(\alpha (1-\alpha) +\frac{(2\alpha_m(1-\alpha_m))^{1/2}}{\pi^{1/2}  } \right)$ (and similarly for $1-\mathbb{P}[B_k]$). Therefore, 
$$\mathbb{P}[|v_{\alpha_m}(L(\widehat{x}_k)) - v_{\alpha_m}(L(x_k^*))|>\frac{\delta}{4}]<\frac{648 (k+|D_0|)}{r_k  \delta^2 (f^*)^2} \left(\alpha (1-\alpha) +\frac{(2\alpha_m(1-\alpha_m))^{1/2}}{\pi^{1/2}  } \right).$$

As a result, 
$$\mathbb{P}[|\mathcal{Y}_m(\widehat{x}_k)-  v_{\alpha_m}(L(x_k^*))|>\frac{\delta}{2}]<
\frac{712 (k+|D_0|)}{r_k  \delta^2 (f^*)^2} \left(\alpha (1-\alpha) +\frac{(2\alpha_m(1-\alpha_m))^{1/2}}{\pi^{1/2}  }\right).$$ By Assumption 4.1, we have $ \sum_{k=1}^{\infty} \frac{k }{r_k }<\infty$. With this assumption, we see that $\sum_{k=1}^{\infty} \frac{|D_0| }{r_k }=|D_0|\sum_{k=1}^{\infty} \frac{1 }{r_k }< \infty$. Therefore,
\[
\begin{split}
\sum_{k=1}^{\infty} \mathbb{P}[|\mathcal{Y}_m(\widehat{x}_k)-  v_{\alpha_m}(L(x_k^*))|>\frac{\delta}{2}] &< \frac{712 }{r_k  \delta^2 (f^*)^2} \left(\alpha (1-\alpha) +\frac{(2\alpha_m(1-\alpha_m))^{1/2}}{\pi^{1/2}  }\right)\sum_{k=1}^{\infty} \frac{(k+|D_0|) }{r_k }   < \infty.
\end{split}
\] 
It follows that $\lim_{n\rightarrow\infty} \mathbb{P}(\cup_{k=n}^\infty \{|\widehat{\mathcal{Y}}_k - v_{\alpha_m}(L(x_k^*))|>\frac{\delta}{2}\})=0$.

\subsection{Proof that $v_{\alpha_m}(L(x_k^*)) \rightarrow  v_{\alpha_m}(L(x^*))$ w.p.1 as $k \rightarrow \infty$} 
According to Theorem 1.3 from \citet{torn1989global}, for a deterministic search (given starting point $x_0$, the design points are determined), the algorithm converges if the design points are everywhere dense, i.e., $v_{\alpha_m}(L(x_k^*))\rightarrow v_{\alpha_m}(L(x^*))$ if $D_k$ is dense in $\mathcal{X}$. When the design points are random ($x_k$ and $x_k^*$ are random variables), denseness of the design points is not sufficient to guarantee almost sure convergence. We next prove that in our algorithm, $v_{\alpha_m}(L(x_k^*))\rightarrow v_{\alpha_m}(L(x^*))$ as $k \rightarrow \infty$ w.p.1. Equivalently, we prove that for all $\delta>0$, $\mathbb{P}[|v_{\alpha_m}(L(x_k^*))- v_{\alpha_m}(L(x^*))|>\delta, i.o.]=0$ ($i.o.$ is shorthand for infinitely often).

For $\epsilon >0$, we can select a region $S$ around $x^*$ such that for all $x\in S$, $|v_{\alpha_m}(L(x))- v_{\alpha_m}(L(x^*))|\leq\epsilon$ (under the assumption that the baseline function $v_{\alpha_m}(L(x))$ is continuous). We next prove that there exists a large value $K_1$ such that at least one design point is selected in $S$ before iteration $K_1$. It then follows that $\mathbb{P}[|v_{\alpha_m}(L(x_k^*))- v_{\alpha_m}(L(x^*))|>\epsilon, i.o.]=0$.

If any points in $S$ are selected in some iteration $k\leq K_1-1$, the condition holds. Now suppose no points in $S$ are selected before iteration $K_1$. In this case we can find a lower bound $\widehat{s}_0^2$ for the predictive variance $\widehat{s}_{K_1}^2(x^*)$ at $x^*$, which is the value of $\widehat{s}_{K_1}^2(x^*)$ if all the design points in $\mathcal{X}\setminus S$ are observed with no noise. Subsequently, we see
$$ T_{K_1}(x^*)> \widehat{s}_0\phi(\frac{-2M}{\widehat{s}_0})-2M\Phi (\frac{-2M}{\widehat{s}_0}) := t_0.$$

In other words, we can find a lower bound for the EI function value at $x^*$, $t_0$. Note that $t_0$ is the EI function value if the predictive response value is $2M$ larger than the current best value and the predictive variance is $\widehat{s}_0$. As the EI function is always positive if $\widehat{s}_0>0$, we see that $t_0>0$. From the proof in Appendix F, we see that if we keep reducing the value of $c$ in the stopping rule to $t_0$, then within a finite number of iterations, the EI function values at all points in $\mathcal{X}\setminus S$ will become smaller than $t_0$. As a result, when $K_1$ is large enough, $T_{K_1}(x)<t_0$ for all $x\in\mathcal{X}\setminus S$ while ${T}_{K_1}(x^*)>t_0$. Therefore, the next design point must belong to $S$, and we finish the proof.

\subsection{Proof that $\widehat{\mathcal{Y}}_k \rightarrow v_{\alpha_m}(L(x^*))$ w.p.1 as $k\rightarrow \infty$}
Since 
$$\{|\widehat{\mathcal{Y}}_k -v_{\alpha_m}(L(x^*))|>\delta \} \subset \{|\widehat{\mathcal{Y}}_k - v_{\alpha_m}(L(x_k^*))|>\frac{\delta}{2}\} \cup \{|v_{\alpha_m}(L(x_k^*))- v_{\alpha_m}(L(x^*))|>\frac{\delta}{2}\},$$
it follows that 
$$\cup_{k=n}^\infty \{|\widehat{\mathcal{Y}}_k -v_{\alpha_m}(L(x^*))|>\delta \} \subset \left\{ \cup_{k=n}^\infty \{|\widehat{\mathcal{Y}}_k - v_{\alpha_m}(L(x_k^*))|>\frac{\delta}{2}\}\right\} \cup  \left\{ \cup_{k=n}^\infty \{|v_{\alpha_m}(L(x_k^*))- v_{\alpha_m}(L(x^*))|>\frac{\delta}{2}\}\right\}.$$

From Appendix G.1, we have that, for all $\delta>0$,
$$\lim_{n\rightarrow\infty} \mathbb{P}(\cup_{k=n}^\infty \{|\widehat{\mathcal{Y}}_k - v_{\alpha_m}(L(x_k^*))|>\frac{\delta}{2}\})=0.$$ Moreover, since $v_{\alpha_m}(L(x_k^*))\rightarrow v_{\alpha_m}(L(x^*))$ as $k \rightarrow \infty$ w.p.1 (Appendix G.2), we have that, for all $\delta>0$,
$$\lim_{n\rightarrow\infty} \mathbb{P}(\cup_{k=n}^\infty \{|v_{\alpha_m}(L(x_k^*))- v_{\alpha_m}(L(x^*))|>\frac{\delta}{2}\})=0 . $$  
Therefore,
\[
\begin{split}
&\lim_{n\rightarrow\infty}\mathbb{P}[\cup_{k=n}^\infty \{|\widehat{\mathcal{Y}}_k -v_{\alpha_m}(L(x^*))|>\delta \} ] \\ \leq &\lim_{n\rightarrow\infty}\mathbb{P}[\{ \cup_{k=n}^\infty \{|\widehat{\mathcal{Y}}_k - v_{\alpha_m}(L(x_k^*))|>\frac{\delta}{2}\}\}] +\lim_{n\rightarrow\infty}\mathbb{P}[  \{ \cup_{k=n}^\infty \{|v_{\alpha_m}(L(x_k^*))- v_{\alpha_m}(L(x^*))|>\frac{\delta}{2}\}\}]=0.
\end{split}
\]
It follows that $\widehat{\mathcal{Y}}_k \rightarrow v_{\alpha_m}(L(x^*))$ w.p.1 as $k\rightarrow \infty$. 
% Appendix here
% Options are (1) APPENDIX (with or without general title) or 
%             (2) APPENDICES (if it has more than one unrelated sections)
% Outcomment the appropriate case if necessary
%
% \begin{APPENDIX}{<Title of the Appendix>}
% \end{APPENDIX}
%
%   or 
%

% CASE 1: BiBTeX used to constantly update the references 
%   (while the paper is being written).
\bibliographystyle{informs2014} % outcomment this and next line in Case 1
\bibliography{sample} % if more than one, comma separated

% CASE 2: BiBTeX used to generate mypaper.bbl (to be further fine tuned)
%\input{mypaper.bbl} % outcomment this line in Case 2

\end{document}